
\input amstex
\documentstyle{amsppt}
\language0
\let\totoc\relax
\redefine\o{\circ}
\define\al{\alpha}
\define\be{\beta}
\define\ga{\gamma}
\define\de{\delta}

\define\ka{\kappa}
\define\la{\lambda}
\define\rh{\rho}
\define\si{\sigma}

\define\ph{\varphi}

\define\ps{\psi}
\define\om{\omega}
\define\Ga{\Gamma}

\define\La{\Lambda}

\define\Om{\Omega}

\define\row#1#2#3{#1_{#2},\ldots,#1_{#3}}
\define\rowup#1#2#3{#1^{#2},\ldots,#1^{#3}}
\define\rwedge#1#2{dx^{#1}\wedge \dots \wedge dx^{#2}}
\define\vvdf#1#2#3{dx^{#1}\wedge \dots \wedge dx^{#2}\otimes \tfrac
\partial{\partial x^{#3}}}
\define\x{\times}
\def\nmb#1#2{#2}            

\define\Mf{\Cal Mf}

\redefine\L{{\Cal L}}

\redefine\lim{\operatorname{lim}}

\define\inv{\operatorname{inv}}

\define\supp{\operatorname{supp}}

\define\dd#1{\tfrac \partial{\partial #1}}
\def\today{\ifcase\month\or
 January\or February\or March\or April\or May\or June\or
 July\or August\or September\or October\or November\or December\fi
 \space\number\day, \number\year}

\topmatter
\title On multilinear operators commuting with Lie derivatives \endtitle
\author
Andreas \v Cap\\
Jan Slov\'ak
\endauthor
\address Institut f\"ur Mathematik, Universit\"at Wien, Strudlhofgasse
4, 1090 Wien, Austria\newline
and Erwin Schr\"odinger International Institute for Mathematical Physics,
Pasteurgasse 6/7, 1090 Wien, Austria\endaddress
\address Department of Algebra and Geometry, Masaryk University,
Jan\'a\v ckovo n\'am.~2a\newline 662~95~Brno, Czech Republic
\endaddress
\abstract
Let $E_1,\dots ,E_k$ and $E$ be natural vector bundles defined over
the category $\Cal Mf_m^+$ of smooth oriented $m$--dimensional
manifolds and orientation preserving local diffeomorphisms, with
$m\geq 2$. Let $M$ be an object of $\Cal Mf_m^+$ which is connected.
We give a complete classification of all separately continuous
$k$--linear operators
$D\:\Ga _c(E_1M)\x\dots\x\Ga _c(E_kM)\to \Ga (EM)$ defined on
sections with compact supports, which commute with
Lie derivatives, i\.e\. which satisfy
$$\Cal L_X(D(s_1,\dots ,s_k))=\sum _{i=1}^kD(s_1,\dots
,\Cal L_Xs_i,\dots ,s_k),$$
for all vector fields $X$ on $M$ and sections $s_j\in\Ga_c(E_jM)$, in
terms of local natural operators and absolutely invariant sections.
In special cases we do not need the continuity assumption. We
also present several applications in concrete geometrical
situations, in particular we give a completely algebraic
characterization of some well known Lie brackets.
\endabstract
\thanks First author supported by project P 77724 PHY of `Fonds zur
F\"orderung der wissenschaftlichen Forschung'. Second author supported
by a grant of the GA\v CR.
\endthanks
\endtopmatter
\document
\head\totoc\nmb0{1}. Natural operators \endhead
In this section we give a brief survey over notions and results
from the theory of natural bundles and operators which we will need
in the sequel and we formulate our main result.
Besides the references cited for the results they can
all be found in \cite{Kol\'a\v r--Michor--Slov\'ak}.

\subheading{\nmb.{1.1}. Definition} A {\sl bundle functor\/} (or
natural bundle) on the category
$\Cal Mf_m^+$ of $m$-dimensional oriented manifolds and
orientation preserving smooth locally invertible mappings (so called
local diffeomorphisms) is a functor $F\:\Cal Mf_m^+\to \Cal F\Cal M$
with values in fibered manifolds which satisfies:

{\narrower
\noindent (i) $B\o F=\text{id}_{\Cal Mf_m^+}$ where
$B\:\Cal F\Cal M\to \Cal Mf$ is the base functor

\noindent (ii) for every inclusion $i:U\to M$ of an open submanifold, $FU$ is
the restriction $p_M^{-1}(U)$ of the value $FM= (p_M\:FM\to M)$ to
$U$ and $Fi$ is the inclusion $p_M^{-1}(U)\to FM$.
\par}

A {\sl vector bundle functor\/} (or natural vector bundle)
is a bundle functor with values in
the category of finite dimensional vector bundles and fiberwise
invertible vector bundle homomorphisms.

\subheading{\nmb.{1.2}} By the general theory, cf\.
\cite{Palais--Terng}, \cite{Epstein--Thurston}, each natural bundle is
obtained (up to isomorphisms) by the construction of associated
bundles to the (higher order) frame bundles $P^{r^+}M$. The latter
principal bundles  are defined as the bundles of $r$-jets of locally
defined orientation preserving diffeomorphisms
$\Bbb R^m\to M$ at the origin $0\in \Bbb R^m$, and their structure
groups are the Lie groups of $r$-jets of orientation preserving
diffeomorphisms
$\Bbb R^m\to \Bbb R^m$ which keep the origin fixed, the so called
{\sl jet groups}. Then the values of $F$ on morphisms depend only on
$r$-jets for some non-negative integer $r$, which is called the {\sl
order of $F$\/}.

In particular, all natural bundles transform smoothly parameterized
families of morphisms into smoothly parameterized families (the so
called regularity). Hence we can apply a bundle functor $F$ to flows of
vector fields on $M$ and we get flows on the values $FM$. This gives
rise to the so called {\sl flow operator \/} which transforms vector
fields on  $M$ into vector fields on $FM$.

Vector bundle functors correspond to linear representations of
the jet groups in finite dimensional vector spaces and the flow
operators can be used to define Lie derivatives of
sections of a natural vector bundle analogously to the classical
case of tensor bundles, cf\. \cite{Terng} and \cite{\v Cap--Slov\'ak}.
By  the regularity of vector bundle functors, the Lie
derivatives define a continuous action of the Lie algebra of
vector fields on the space of sections. This action is local in both
arguments and so in view of the localization property of the vector
bundle functors we get for each vector bundle functor $F$
an action of the sheaf of vector
fields on an oriented manifold $M$ on the sheaf of local sections of
$FM$.

\subheading{\nmb.{1.3}. Absolutely invariant sections}
Let $F$ be a natural vector bundle defined on the category
$\Cal Mf_m^+$ and let $M$ be an object of $\Cal Mf_m^+$ which is
connected. A section $s\in \Ga (FM)$ is called {\sl absolutely
invariant\/} if and only if $\Cal L_Xs=0$ for all vector fields $X$
on $M$. The finiteness of the order of the bundle functor $F$
and arguments from the proof of the main theorem of \cite{\v Cap--Slov\'ak}
imply that $s$ is then invariant for the actions of all orientation
preserving diffeomorphisms of $M$. Now assume that $F$ is of order $r$ and
let $\rho :G^{r+}_m\to GL(V)$ be the corresponding representation of
the jet group. Then one easily verifies that the absolutely invariant
sections of $FM$ are in bijective correspondence with elements of $V$
which are invariant for the action $\rho$. In particular we see that
the absolutely invariant sections always form a finite dimensional
vector space of dimension at most the fiber dimension of $FM$, and
that an absolutely invariant section is determined by its restriction
to one point. Moreover an absolutely invariant section $s$ of $FM$
gives rise to a unique absolutely invariant section of $FN$ for any
oriented $m$--dimensional manifold $N$.

Note that in many concrete situations there exist no absolutely
invariant sections. For example in a subbundle of a tensor bundle
$\otimes^pTM\otimes\otimes^qT^*M$ absolutely invariant sections can
exist only if $p=q$.

\subheading{\nmb.{1.4}. Definition} Let $E,
F\:\Cal Mf_m^+\to \Cal F\Cal M$ be two bundle functors. A {\sl natural
operator\/} $D\:E\to F$ is a system $D_M$ of local smooth operators, for all
oriented $m$-dimensional manifolds $M$, which commute with the
actions of local diffeomorphisms. So $D_M\:\Ga(EM)\to\Ga(FM)$
transforms smooth families of sections of $FM$ into smooth families
of sections of $GM$ and $D_M(s)(x)$ depends only on the germ of $s$ at $x$
for all sections $s\in \Ga(EM)$ and points $x\in M$. Moreover for any
diffeomorphism $\ph :M\to N$ we have
$F(\ph )_*\o D_M=D_N\o E(\ph )_*\:\Ga(EM)\to \Ga(FN)$, and for any
open submanifold $U\subset M$ and section $s\in\Ga (EM)$ we have
$D_U(s|_U)=D_M(s)|_U$.

Each local $k$-linear natural operator $D\:E_1\x \dots\x
E_k\to E$ satisfies
$$\Cal L_X(D_M(s_1,\dots ,s_k))=
\sum _{i=1}^kD_M(s_1,\dots ,\Cal L_Xs_i,\dots ,s_k)$$
for each manifold $M$, all vector fields $X$ on $M$ and all
$s_i\in \Ga(E_iM)$. We say that $D$ {\sl commutes with Lie derivatives}.
Moreover, each $k$-linear local operator $D_M\:\Ga(E_1M)\x \dots\x
\Ga(E_kM)\to \Ga(EM)$ which commutes with Lie derivatives extends
to a unique natural operator $D$, see \cite{\v Cap--Slov\'ak, 3.2}.

\subheading{\nmb.{1.5}. Determination of local natural operators}
Let us consider a $k$--linear natural operator
$D:E_1\x \dots\x E_k\to E$, where the $E_i$ and $E$ are natural
vector bundles. By the multilinear Peetre theorem (c\.f\.
\cite{Cahen--de~Wilde--Gutt} or \cite{Slov\'ak}) $D$
is of some finite order $r$. Let $J^rE_i$ denote the
$r$-jet prolongations of the natural vector bundles $E_i$, let $V_i$
and $V$ be the standard fibers of the $E_i$ and $E$, and let
$T^r_mV_i$ be the standard fiber of $J^rE_i$. The natural operator $D$
is uniquely determined by its value on $\Bbb R^m$ and $D_{\Bbb R^m}$ is
induced by a $k$--linear vector bundle homomorphism
$J^r E_1\Bbb R^m\x\dots\x J^r E_k\Bbb R^m\to E\Bbb R^m$. This
in turn is uniquely determined by the induced linear map $\tilde D\:
T^r_mV_1\otimes\dots\otimes T^r_mV_k\to V$. The spaces $T^r_mV_i$ and
$V$ are modules over the jet group $G^{r+\ell}_m$, where $\ell$ is
the maximum of the orders of the bundle functors $E_i$ and $E$, and
the fact that $\tilde D$ comes from a natural operator is equivalent
to equivariancy of $\tilde D$ for these actions.

Thus the determination of natural operators of fixed order $r$
reduces to the classification of equivariant linear maps between two
finite dimensional representations of the jet group $G^{r+\ell}_m$.
This classification has been carried out in a variety of special cases,
see \cite{Kol\'a\v r--Michor--Slov\'ak} for a collection of results in this
direction.

In particular the space of multilinear local natural operators of
fixed order between fixed bundles is always finite dimensional.
Moreover we can prove that there is always some maximal order for
such operators:

\proclaim{\nmb.{1.6}. Lemma }
Let $E_1,\dots ,E_k$ and $E$ be vector bundle functors on
$\Cal Mf_m^+$. Then the space of all $k$--linear (local) natural
operators $D:E_1\x\dots\x E_k\to E$ is finite dimensional.
\endproclaim
\demo{Proof}
As discussed above, each such operator
$D_M\:\Ga(E_1M)\oplus\dots\oplus \Ga(E_kM)\to \Ga(EM)$ factors
through a $k$--linear $G^{r+\ell}_m$--equivariant map on the jets
of sections in one point of $\Bbb R^m$, $\tilde D\:
T^r_mV_1\otimes\dots\otimes T^r_mV_k\to V$, where $r$ is the order of
the operator, $\ell$ is the maximum of the orders of the natural
bundles. Let us write $D^*\:
V^*\to (T^r_mV_1)^*\otimes\dots\otimes (T^r_mV_k)^*$ for the
dual linear mapping to the corresponding linear mapping  on the
tensor product of the jet spaces. The jet group $G^1_m=GL(m,\Bbb R)$
of order one is always a subgroup and thus, in particular the mapping
$D^*$ is an $SL(m,\Bbb R)$--module homomorphism. As representation
spaces for $SL(m,\Bbb R)$, all the standard fibers $V_1,\dots, V_k$
and $V$ decompose into irreducible representations, as well as the
whole tensor product of the duals. Each of the irreducible components
is determined by a highest weight vector sitting in the whole space,
we shall write $V_\la$ for a component with highest weight $\la$.
Furthermore, we always have $T^r_m W=\oplus_{i=0}^r
W\otimes S^i\Bbb R^{m*}$ as $GL(m,\Bbb R)$--modules, in particular
$T^r_m(W_1\oplus W_2)=T^r_mW_1 \oplus T^r_mW_2$ for any two
$GL(m,\Bbb R)$--modules. Thus we may rewrite the target space of $D^*$ in
the form
$$\align
\bigotimes_{i=1}^k\biggl( \bigoplus_{j=0}^r
V_{i}^*\otimes S^j \Bbb R^{m} \biggr) &= \bigoplus_{s=1}^N
\bigoplus_{j_1,\dots,j_k=0}^r\biggl( V_{\rh_s}\otimes\otimes_{i=1}^k
S^{j_i}\Bbb R^{m}\biggr)\\
&\subset
\bigoplus_{s=1}^N
\bigoplus_{j_1,\dots,j_k=0}^r\biggl(
V_{\rh_s}\otimes\otimes^{j_1+\dots+j_k}\Bbb R^{m}\biggr)
\endalign$$
where the spaces $V_{\rh_s}$ are the irreducible components in the
tensor product $V_1^*\otimes\dots\otimes V_k^*$. Since the decomposition
into irreducible components is always unique up to order and since a
homomorphism between irreducible components is a multiple of the
identity, the mapping $D^*$ is just the inclusion of the irreducible
components of $V^*= \oplus_i V_{\mu_i}$ into the tensor product. Thus
we have only to show that highest weight vectors in the tensor product
with the weights $\mu_i$ can appear only in the components
$V_{\rh_s}\otimes\otimes^{j}\Bbb R^{m}$ with $j$ less that some fixed
bound.

If our jet spaces were completely reducible $GL(m,\Bbb R)$--modules,
we could apply just the above discussion and to compare the
actions of the center of $GL(m,\Bbb R)$ on the domain and the target.
The multiples of the identity matrix in the Lie algebra
$\goth g\goth l(m,\Bbb R)$ must act by some scalar and each
copy of $\Bbb R^m$ in the tensor product increases this scalar  by
one. This simple idea leads then directly to a bound as required.
Since we do not want to restrict ourselves to completely decomposable
$GL(m,\Bbb R)$--modules, we have to use the above decomposition into
irreducible $SL(m,\Bbb R)$--modules and to apply the equivariancy
with respect to the induced action of the universal enveloping
algebra of $\goth s\goth l(m,\Bbb R)$, where the center is not
trivial.

Let us assume there is a highest weight vector with weight $\mu$ in
$V_{\rh}\otimes\otimes^{j}\Bbb R^{m}$. Then we can complexify both
the representation spaces, and we shall find a highest weight vector
with weight $\mu$ in the $\goth U(\goth s\goth l(m,\Bbb C))$--module
$V_{\rh_s}\otimes\otimes^{j}\Bbb C^{m}$ (now $V_{\rh_s}$ means the
complex irreducible representation with highest weight $\rh_s$).

For the semisimple Lie algebra $\goth s\goth l(m,\Bbb C)$, there is the
Casimir operator $C$
sitting in the center of $\goth U(\goth s\goth l(m,\Bbb C))$ and so
acting by a scalar $C_\mu$ on each irreducible representation space
$V_\mu$. This scalar can be computed easily by  means of the highest
weights,
$C_\mu = \langle \mu, \mu+2\de\rangle$ where $2\de$ is the sum of all
positive roots  and $\langle \ ,\ \rangle$ is the Killing form on the
algebra, see e.g\. \cite{Samelson, p\. 121}. In our case,
the weights are in the dual to the Cartan algebra which consists of
trace free diagonal matrices. If $e^1,\dots, e^{m}$ is the usual
linear basis, the positive roots are $e^i-e^j$ with $i<j$, all
possible highest weights are of the form
$a_1e^1+\dots+a_{m-1}e^{m-1}$ with integers
$a_1\ge a_2\ge \dots\ge a_{m-1}\ge 0$, and the Killing form is just the
Euclidean metric up to a constant negative multiple. In particular we
get $\de=(m-1)e^1 + (m-2)e^2 + \dots + e^{m-1}$.

In order to estimate
$C_\mu$ for our irreducible component in the tensor product, we have
to know which weights may appear in a tensor product of two
irreducible representations. A
general answer is obtained easily
by means of some of the well known consequences of the
Weyl's character formula, e.g\. Klimyk's formula, see
\cite{Samelson, p\. 128}: We have only to consider the highest
weight of the first representation and to add any weight of the other
one (only those which yield dominant weights can apply and the
multiplicity must be computed). In our case this means, we take the
weight $\rh$ and we add, step by step, the weights of $\Bbb C^m$. To
get all of them, we apply first the Weyl group  to the highest weight
$e^1$ of $\Bbb C^m$ which yields just the weights $e^1,\dots,e^m$.
These must be involved, but the dimension of the sum of the
corresponding weight spaces is at least $m$, so we are done.
Altogether, a highest weight in $V_{\rh}\otimes\otimes^{j}\Bbb R^{m}$
must be of the form $\rh + e^{i_1}+ \dots+ e^{i_j}$ and we obtain
$$\align
C_\mu &= \langle  \rh + e^{i_1}+ \dots+ e^{i_j}, \rh +
e^{i_1}+ \dots+ e^{i_j} + 2((m-1)e^1 + \dots + e^{m-1})\rangle \\
&=C_\rh + 2\langle  \rh, e^{i_1}+ \dots+ e^{i_j}\rangle +
\langle e^{i_1}+ \dots+ e^{i_j},e^{i_1}+ \dots+ e^{i_j}+2\de\rangle \\
&\le \text{const} - \text{a positive constant multiple of $j$.}
\endalign$$
This estimate implies that for each fixed $\rh$ and $\mu$, the order
$j$ is bounded and the lemma is proved.
\qed\enddemo

In fact the proof of lemma \nmb!{1.6} shows even that the same
statement holds for vector bundle functors and operators defined on
the category of smooth oriented $m$--dimensional manifolds with fixed
volume forms and volume preserving local diffeomorphisms.

Next we introduce a class of non--local operators which commute with
Lie derivatives.

\subheading{\nmb.{1.7}. Definition }
Let $E_1,\dots ,E_k,E$ be natural vector bundles defined on the
category $\Mf_m^+$, and let us write $\Ga_c(E_iM)$ for the space of
compactly supported smooth sections of the vector bundle $E_iM$ over
an $m$-dimensional oriented manifold $M$ while
$\Ga(EM)$ means the space of all smooth sections. Let $I^1=\{
i^1_1<\dots <i^1_{n_1}\},\dots ,I^r=\{ i^r_1<\dots <i^r_{n_r}\},
J=\{ j_1<\dots <j_n\}$ with $i^1_1<\dots <i^r_1$ be a partition of
the set $\{ 1,\dots ,k\}$ into disjoint subsets, such that all $I^j$
are nonempty, while $J$ is allowed to be empty.

An {\sl elementary almost natural operator\/}
$\Ga_c(E_1M)\x\dots \x\Ga_c(E_kM)\to \Ga(EM)$
{\sl of type\/} $(I^1,\dots ,I^r,J)$ is an operator of the form
$$(s_1,\dots s_k)\mapsto \la^1(s_{i^1_1},\dots ,s_{i^1_{n_1}})\dots
\la^r(s_{i^r_1},\dots ,s_{i^r_{n_r}})D_M(s_{j_1},\dots ,s_{j_n}),$$
where $D:E_{j_1}\x\dots \x E_{j_n}\to E$ is a $j_n$--linear local
natural operator and any $\la^\ell$ is an operator
$\Ga_c(E_{i^\ell_1})\x\dots\x\Ga_c(E_{i^\ell_{n_\ell}})\to \Bbb R$
of the form
$$(s_1,\dots ,s_{n_\ell})\mapsto \int_M\langle D^\ell_M(s_1,\dots ,
s_{n_{\ell}-1}),s_{n_\ell}\rangle,$$
where $D^\ell:E_{i^\ell_1}\x\dots\x E_{i^\ell_{n_\ell-1}}\to
E^*_{i^\ell_{n_\ell}}\otimes \La^mT^*$ is a $(n_\ell-1)$--linear
local natural operator
and $\langle \quad ,\quad \rangle$ denotes the canonical pairing
$\Ga_c(E^*_{i^\ell_{n_\ell}}M\otimes \La^mT^*M)\x
\Ga_c(E_{i^\ell_{n_\ell}}M)\to \Ga_c(\La^mT^*M)$.
In this notation we use the convention that a local natural operator
without arguments is just an absolutely invariant section of the
target bundle.

An {\sl almost natural operator\/} $E_1M\x\dots \x E_kM\to EM$ is a
linear combination of elementary almost natural operators.

\proclaim{\nmb.{1.8}. Proposition }
Any almost natural operator commutes with Lie derivatives.
\endproclaim
\demo{Proof}
Let us first show that for a multilinear map $\la$ as in definition
\nmb!{1.7} we have $\sum_{i=1}^n\la(s_1,\dots,\L_Xs_i,\dots,s_n)=0$
for any sections $s_i$ and any vector field $X$ on $M$.
Since any Lie derivative of an $m$-form on an $m$-dimensional
manifold is an exact form we have:
$$\align
0&=\int_M\L_X(\langle D_M(s_1,\dots ,s_{n-1}),s_n\rangle)=\\
&=\int_M(\langle \L_X(D_M(s_1,\dots ,s_{n-1})),s_n\rangle+\langle D_M(s_1,\dots
,s_{n-1}),\L_Xs_n\rangle )=\\
&=\sum_{i=1}^{n-1}\int_M\langle D_M(s_1,\dots,\L_Xs_i,\dots
,s_{n-1}),s_n\rangle +\\
&\qquad\qquad\qquad\qquad+\int_M\langle D_M(s_1,\dots ,s_{n-1}),\L_Xs_n\rangle
=\\
&=\sum_{i=1}^n\la (s_1,\dots ,\L_Xs_i,\dots ,s_n)
\endalign$$
{}From this the result is obvious.
\qed\enddemo

The aim of this paper is to discuss to what extent the converse of
this proposition holds, i.e\. whether all multilinear operators which
commute with Lie derivatives are almost natural.

The space of sections of a vector bundle has a canonical Fr\'echet
topology (in principle just uniform convergence in each derivative),
while the space of compactly supported sections of a bundle can be
canonically topologized as the inductive limit of the spaces of
sections having support in some fixed compact subset, which are again
Fr\'echet spaces (see \nmb!{2.2} for a more detailed description).
Obviously any almost natural operator is jointly continuous for these
topologies, i.e\. continuous as a mapping $D:\Ga_c(E_1M)\x\dots
\x\Ga_c(E_kM)\to \Ga(EM)$.

The main result of this paper is that under the weaker assumption of
continuity in each argument separately, the converse to \nmb!{1.8} is
actually true:
\proclaim{\nmb.{1.9}. Theorem }
Let $E_1,\dots ,E_k,E$ be natural vector bundles defined on the
category $\Mf_m^+$, $m\geq 2$, and let $M$ be a connected oriented
smooth manifold of dimension $m$. Then any separately continuous
$k$--linear operator $D:\Ga_c(E_1M)\x\dots \x\Ga_c(E_kM)\to \Ga(EM)$
which commutes with Lie derivatives is an almost natural operator. In
particular any such operator is automatically jointly continuous and
the space of such operators is always finite dimensional and independent
of the manifold $M$.
\endproclaim

Let us remark that the assumption of separate continuity is needed
in one step in the proof, the application of the Schwartz kernel
theorem in \nmb!{3.3}. In our opinion, the only way to avoid the
continuity assumption is to prove that continuity already follows
from the commutation with Lie derivatives. Certainly this is a rather
difficult problem, since there even are examples of (non--linear)
local natural operators which are not continuous, see e.g\.
\cite{Kol\'a\v r--Michor--Slov\'ak, 23.9}.

\head\totoc\nmb0{2}. Natural presheaves \endhead
Although our results apply mainly to natural vector
bundles and their sections we have to use the more
technical notion of an admissible natural presheaf in the proofs.

\subheading{\nmb.{2.1}. Definitions}
Let $M$ be a smooth manifold. By $\Cal X(\quad )$ we denote the sheaf
of smooth vector fields over $M$.

\noindent (1) A {\sl natural presheaf\/} $\Cal F$ on $M$ is a presheaf of
locally convex modules over the sheaf of Lie algebras $\Cal X(M)$. This means
that for
any open subset $U$ of $M$ there is a locally convex vector space
$\Cal F(U)$ on which the Lie algebra $\Cal X(U)$ acts by continuous
linear operators, which we denote by $\Cal L_X$ for $X\in \Cal X(U)$.
By the same symbol we denote the actions of the restrictions of $X$ to
open subsets of $U$ on the corresponding spaces.
Moreover for any open subset $V\subset U$ there is a continuous
linear restriction mapping $r^U_V:\Cal F(U)\to \Cal F(V)$ which is
equivariant over the corresponding restriction map for vector fields
for the Lie algebra actions.

This is an infinitesimal version of the definition of a natural
presheaf given in \cite{Eck}.

\noindent (2) An element $s\in \Cal F(U)$ is called {\sl absolutely
invariant\/} if and only if $\Cal L_Xs=0$ for all $X\in \Cal X(M)$
with support contained in $U$.

\noindent(3) A natural presheaf $\Cal F$ over $M$ is called
{\sl admissible\/} if and only if the following two conditions are
satisfied:

{\narrower
\noindent(i) For any open subset $U$ of $M$ and system of
countably many open sets $\{ U_i\}$ such that $U_i\subset U_{i+1}$
for all $i$, and such that $U=\cup _{i\in \Bbb N}U_i$,
each element of $\Cal F(U)$ is uniquely determined by
its restrictions to all spaces $\Cal F(U_i)$.

\noindent(ii) For any connected open subset $U$
of $M$ and any absolutely invariant element $s\in \Cal F(U)$ there is a
unique absolutely invariant element $\tilde s\in \Cal F(M)$ such that
$r^M_U(\tilde s)=s$. Moreover the operator of extending absolutely
invariant elements is continuous.

}
\noindent(4) Let $E_1\to M,\dots ,E_k\to M$ be vector bundles, $\Cal F$ a
presheaf of vector spaces over $M$ and let us write $\Ga _c(E_i)$
for the space of sections of $E_i$ with compact support. A $k$--linear operator
$D:\Ga _c(E_1)\x\dots\x \Ga _c(E_k)\to \Cal F(M)$ is
called {\sl local\/} if and only if for each open subset $U$ of $M$
and each $i$, $1\leq i\leq k$, the condition $s_i|_U=0$ implies that
$r^M_U(D(s_1,\dots ,s_k))=0$.

\noindent(5) Let $E_1\to M,\dots ,E_k\to M$ be vector bundles, $V$ a
vector space. A $k$--linear map
$\la\:\Ga _c(E_1)\x\dots\x \Ga _c(E_k)\to V$ is called {\sl weakly
local\/} if and only if $\la (s_1,\dots ,s_k)=0$ for all sections with
$\supp (s_1)\cap \dots \cap\supp (s_k)=\emptyset$.

\subheading{\nmb.{2.2}. Topologies on spaces of sections and
operators between them}

\noindent Before we give examples of natural presheaves we describe the
topologies we will use. For an open subset $U\subset \Bbb R^m$ we put
on $C^\infty(U,\Bbb R^n)$ as usual the topology of uniform convergence
in all derivatives separately, which is a Fr\'echet topology. Next for
a smooth vector bundle $E\to M$ over a smooth manifold a vector
bundle atlas $(U_i)$ gives rise to linear maps
$\Ga (E)\to C^\infty(U_i,\Bbb R^n)$, where $n$ is the fiber dimension
of $E$, and we put on $\Ga (E)$ the initial topology with respect to
these maps. Since we can use a countable atlas $\Ga (E)$ is again a
Fr\'echet space. Next for any compact subset $K\subset M$ the space
$\Ga_c(E;K)$ of all sections of $E$ with support contained in $K$ is
a closed linear subspace of $\Ga (E)$ and thus also a Fr\'echet space
with the induced topology. Now for $i\in \Bbb N$ we can choose compact
subsets $K_i\subset M$ such that $K_i\subset K_{i+1}$ and such that
$M=\cup_{i\in \Bbb N}K_i$. Then $\Ga _c(E;K_i)$ is a closed linear
subspace of $\Ga _c(E;K_{i+1})$ and
$\Ga _c(E)=\cup _{i\in \Bbb N}\Ga _c(E;K_i)$ as a set. We
topologize $\Ga _c(E)$ as the inductive limit of the spaces
$\Ga _c(E;K_i)$, so $\Ga _c(E)$ is the strict inductive limit of a
sequence of Fr\'echet spaces. In particular this implies that
$\Ga _c(E)$ is a bornological space, i\.e\. a linear map from
$\Ga _c(E)$ to any locally convex vector space is continuous if and
only if it is bounded.

Consider the space $L(\Ga _c(E_1),\dots ,\Ga _c(E_k);\Ga (E))$ of all
bounded $k$--linear maps
$\Ga _c(E_1)\x \dots\x \Ga _c(E_k)\to \Ga (E)$, where the $E_i$ and
$E$ are vector bundles over $M$.
Since $\Ga _c(E_i)$ is bornological we conclude from \cite{Fr\"olicher--Kriegl,
3.7.5}
that the bounded $k$--linear maps are exactly the separately
continuous ones.

Now we declare a subset
$B\subset L(\Ga _c(E_1),\dots ,\Ga _c(E_k);\Ga (E))$ to be bounded if
and only if $B(A_1\x\dots\x A_k)$ is bounded in $\Ga (E)$ for all
bounded subsets $A_i\subset \Ga _c(E_i)$. Then on
$L(\Ga _c(E_1),\dots ,\Ga _c(E_k);\Ga (E))$ we put the associated
locally convex topology, i\.e\. we take as a basis of neighborhoods
of zero all absolutely convex subsets $U$ such that for any bounded
subset $B$ there exists some $t\in \Bbb R$ such that
$B\subset t\cdot U$. Note that by definition this is a bornological
topology. Moreover by \cite{Fr\"olicher--Kriegl, 3.7.3} flipping coordinates
gives an
isomorphism
$$\multline
L(\Ga _c(E_1),\dots ,\Ga _c(E_k);\Ga (E))\cong\\
\cong L(\Ga _c(E_1),\dots
,\Ga _c(E_\ell);L(\Ga _c(E_{\ell+1}),\dots ,\Ga _c(E_k);\Ga (E))).
\endmultline$$

\subheading{\nmb.{2.3}. Examples of natural presheaves}
(1) Let $E$ be a vector bundle functor defined on the category
$\Cal Mf_m$ of $m$--dimensional smooth manifolds and local
diffeomorphisms. Then with the usual notion of the Lie derivative
along vector fields and the topology described above, the spaces
$\Ga (EU)$
of smooth sections of $EU$ form a natural presheaf (in fact a
natural sheaf) over $M$. In section \nmb!{4} we will
prove that this natural presheaf is admissible.

\noindent(2) A less simple example which will be crucial in the sequel
is the following: Let $E_1,\dots ,E_k,E$  be vector bundle functors as
above. Put
$$\Cal F(U):=L(\Ga _c(E_1U),\dots ,\Ga _c(E_kU);\Ga (EU)),$$
the space of all separately continuous $k$--linear operators.
On this space we define the action of the Lie algebra $\Cal X(U)$ by
$$(\Cal L_XD)(s_1,\dots ,s_k):=\Cal L_X(D(s_1,\dots ,s_k))-\sum_{i=1}^k
D(s_1,\dots ,\Cal L_Xs_i,\dots ,s_k).$$
This action is continuous since it is just given by compositions with
continuous linear operators. To define the restriction
mappings just note that for open sets $V\subset U$ any section with
compact support of the restricted vector bundle over $V$ can be
extended by zero to a section over $U$. It will be shown in section
\nmb!{5} that this natural presheaf is also admissible. Note that
the absolutely invariant elements in $\Cal F(U)$ are by definition
exactly those operators which commute with Lie derivatives.

\subheading{\nmb.{2.4}. Remark }
When dealing with operators with values in spaces of operators one
has to be very careful about the various meanings of locality and
weak locality. Moreover the name `weakly local' is a little
misleading when dealing with operators with values in general
admissible natural presheaves. Clearly any local operator with values
in a sheaf over $M$ is a weakly local mapping. But if one deals with
presheaves which only satisfy the very weak condition (i) of
\nmb!{2.1}.(3) then there are local operators which are
not weakly local. Let
us discuss this in an example which is relevant in the sequel:

For a
smooth manifold $M$ of dimension $m$ consider the trilinear operator
$$D:C^\infty_c(M,\Bbb R)\x C^\infty_c(M,\Bbb R)\x \Om^m_c(M)\to
C^\infty(M,\Bbb R)$$
defined by $D(f,g,\om):=f\cdot\int_Mg\cdot\om$. As a trilinear
operator this is neither local nor weakly local. But we can consider
$D$ in three ways as a linear operator with values in a space of
bilinear operators. In fact all these associated linear operators are
local, since this just means that if one of the three sections $f,g$
and $\om$ vanishes on an open subset $U$ and the other two sections
have support contained in $U$ then $D(f,g,\om)$ vanishes on $U$ and
this is clearly satisfied.

Next we can consider $D$ in three ways as a bilinear operator with
values in a space of linear operators. In particular consider the
operator $(f,g)\mapsto(\om \mapsto D(f,g,\om))$. Then this operator
is local since if either $f$ or $g$ vanish on an open subset $U$ of
$M$ and $\om$ has support contained in $U$ then $D(f,g,\om)$ vanishes
on $U$. On the other hand this operator is not weakly local since $f$
and $g$ having disjoint supports does not imply that $D(f,g,\om)=0$.
Similarly one checks that the operator
$(f,\om)\mapsto(g \mapsto D(f,g,\om))$ is local but not weakly local
while on the other hand $(g,\om)\mapsto(f \mapsto D(f,g,\om))$ is
weakly local but not local.

\subheading{\nmb.{2.5} }
Let $E_1,\dots ,E_k$ be vector bundle functors defined on the category
$\Cal Mf_m^+$, $m\geq 2$, and let $\Cal F$ be an admissible natural
presheaf over an oriented connected $m$--dimensional manifold $M$. Let
$D:\Ga _c(E_1M)\x \dots \x \Ga _c(E_kM)\to \Cal F(M)$ be a $k$--linear operator
which commutes with Lie derivatives, i\.e\. such that
$$\Cal L_X(D(s_1,\dots ,s_k))=
\sum _{i=1}^kD(s_1,\dots ,\Cal L_Xs_i,\dots ,s_k)$$
for all $X\in \Cal X(M)$, $s_i\in \Ga _c(E_iM)$.

First consider sections $s_i\in \Ga _c(E_iM)$ and suppose that over
some open set $U$ all $s_i$ are identically zero. Then $\Cal L_Xs_i$
is globally zero for all vector
fields $X$  with support in $U$ and so the
commutation with Lie derivatives implies that the restriction of
$D(s_1,\dots ,s_k)$ to $U$ is absolutely invariant.

\proclaim{Lemma }
Suppose that $W:=M\setminus (\cup_i \supp(s_i))\neq \emptyset$. Then
there is a unique absolutely invariant element
$D(s_1,\dots ,s_k)_{\inv}\in \Cal F(M)$ such that
$r^M_W(D(s_1,\dots ,s_k)_{\inv})=r^M_W(D(s_1,\dots ,s_k))$.
\endproclaim
\demo{Proof} Since $\Cal F$ is admissible, the absolutely invariant
restrictions of $D(s_1,\dots,s_k)$ to the individual connected
components of $W$ extend uniquely to globally defined absolutely
invariant elements. Thus, we only
have to show that for any two connected components $U$ and $V$ of
$W$ the two absolutely invariant elements in $\Cal F(M)$ determined by
the restrictions of $D(s_1,\dots ,s_k)$ to $U$ and $V$ coincide.
Since $M$ is connected and of dimension at least 2 we can choose
$k+1$ smooth curves $c_1,\dots ,c_{k+1}:[0,1]\to M$, which all start
in $U$, end in $V$ and are pairwise disjoint. Then there are open
sets $U_i$ and $V_i$ for $1\leq i\leq k+1$ such that
$c_i([0,1])\subset U_i\subset \bar U_i\subset V_i$ and such that the
sets $V_i$ are pairwise disjoint. Next for each $i=1,\dots ,k$ choose a
smooth function $f_i\in C^\infty(M,[0,1])$ which is identically 1 on
$\bar U_i$ and vanishes outside of $V_i$ and put
$f_{k+1}:=1-\sum _{i=1}^kf_i$.

By $k$--linearity
$D(s_1,\dots ,s_k)=\sum _ID(f_{i_1}s_1,\dots ,f_{i_k}s_k)$, where $I$
runs over the set of all multiindices $(i_1,\dots ,i_k)$ with
$1\leq i_j\leq k+1$. But now observe that for any such multiindex there
is at least one integer $\ell$, $1\leq \ell \leq k+1$, which does not appear
in the multiindex. Thus all sections $f_{i_1}s_1,\dots ,f_{i_k}s_k$
vanish locally around the curve $c_\ell$ and hence on a connected
open set which intersects both $U$ and $V$. So for each of the
summands the two extensions coincide and hence the same is true for
the sum.
\qed\enddemo

\subheading{\nmb.{2.6} }
For $i=1,\dots ,k$ choose nonempty open subsets $U_i$ and $V_i$ of $M$
such that $U_i\subset \bar U_i\subset V_i$
and such that the sets $\bar V_i$ are pairwise disjoint.
Since $M$ is connected the open subset
$U_{k+1}:=M\setminus (\cup _i\bar V_i)$ is nonempty.
Clearly it is possible to find such a configuration
with all $V_i$ contained in an arbitrary small open subset of $M$.
Then choose smooth functions $f_i\in C^\infty(M,[0,1])$ such that
$f_i$ is identically 1 on $\bar U_i$ and identically zero on
$M\setminus V_i$ and put $f_{k+1}:=1-\sum _if_i$.

Now fix arbitrary sections $s_i\in \Ga _c(E_iM)$. For each multiindex
$I=(i_1,\dots ,i_k)$ of integers $i_1,\dots ,i_k$ between 1 and
$k+1$ there is an open subset of $M$ on which all the sections
$f_{i_1}s_1,\dots,f_{i_k}s_k$ vanish by the choice of the functions.
So by the previous lemma we get an absolutely invariant element
$D(f_{i_1}s_1,\dots,f_{i_k}s_k)_{\inv}\in \Cal F(M)$.
Let us define
$$
\la (s_1,\dots ,s_k)=\la^{f_1,\dots,f_k}(s_1,\dots, s_k) = \sum_I
D(f_{i_1}s_1,\dots,f_{i_k}s_k)_{\inv}.$$
This operator is obviously $k$--linear,
and since the presheaf $\Cal F$ is admissible it
is continuous provided that $D$ is continuous.

\proclaim{\nmb.{2.7}. Lemma}
The operator $\la$ does not depend on the choices of $U_i$, $V_i$ and $f_i$.
Moreover if for some open subset $U$ of $M$ the sections
$s_j\in \Ga_c(E_jM)$ all restrict to zero on $U$, then
$r^M_U(\la(s_1,\dots ,s_k))=r^M_U(D(s_1,\dots ,s_k))$.
\endproclaim
\demo{Proof}
We show that we may replace $(U_j,V_j,f_j)$ by $(U'_j,V'_j,f'_j)$ for
one fixed $j$ if the new set of triples still satisfies the
conditions of \nmb!{2.6}. Clearly, it suffices to do this for $j=1$.
So let us take a second choice $U'_i$, $V'_i$ and $f'_i$ where the data
do not change for $i>1$.

By the multilinearity of $D$, we get
$$\multline
D(f_{i_1}s_1,\dots,f_{i_k}s_k)-D(f'_{i_1}s_1,\dots,f'_{i_k}s_k)=\\
=\sum_jD(f'_{i_1}s_1,\dots,f'_{i_{j-1}}s_{j-1},(f_{i_j}-f'_{i_j})s_j,
f_{i_{j+1}}s_{j+1},\dots,f_{i_k}s_k).
\endmultline$$
By the assumption on $\bar V'_1$ all sections occurring as arguments in this
equation vanish simultaneously on some open subset in $M$,
so we get the analogous equation for the invariant elements
$D(f_{i_1}s_1,\dots,f_{i_k}s_k)_{\inv}-
D(f'_{i_1}s_1,\dots,f'_{i_k}s_k)_{\inv}$. Thus we obtain
$$\multline
\la ^{f_1,\dots,f_k}(s_1,\dots,s_k)-\la ^{f'_1,\dots,f'_k}(s_1,\dots,s_k)=\\
=\sum_{j=1}^k\sum_ID(f'_{i_1}s_1,\dots,f'_{i_{j-1}}s_{j-1},
(f_{i_j}-f'_{i_j})s_j,f_{i_{j+1}}s_{j+1},\dots,f_{i_k}s_k)_{\inv}.
\endmultline$$
But in fact only the terms with either $i_j=1$ or $i_j=k+1$ can
contribute, since otherwise $f_{i_j}=f'_{i_j}$ and so one argument
vanishes globally on $M$. Thus
$$\multline
\la ^{f_1,\dots,f_k}(s_1,\dots,s_k)-\la ^{f'_1,\dots,f'_k}(s_1,\dots,s_k)=
\\
=\sum_{j=1}^k\sum_{(i_1,\dots,i_{j-1},i_{j+1},\dots,i_k)}
\bigl(D(f'_{i_1}s_1,\dots,f'_{i_{j-1}}s_{j-1},(f_{1}-f'_{1})s_j,
f_{i_{j+1}}s_{j+1},\dots,f_{i_k}s_k)_{\inv}
\\
+D(f'_{i_1}s_1,\dots,f'_{i_{j-1}}s_{j-1},(f_{k+1}-f'_{k+1})s_j,
f_{i_{j+1}}s_{j+1},\dots,f_{i_k}s_k)_{\inv}\bigr).
\endmultline$$
Since $f_{k+1}-f'_{k+1}=f'_1-f_1$ by definition, the two terms in the
most inner sum annihilate each other and so the whole expression vanishes.

Finally if all $s_j$ restrict to zero on an open subset $U$ of $M$
then by the first part of this proof we may compute $\la$ using
triples $(U_i,V_i,f_i)$ such that $f_1,\dots ,f_k$ have support
contained in $U$. Then by definition
$\la (s_1,\dots ,s_k)=D(f_{k+1}s_1,\dots ,f_{k+1}s_k)_{\inv}=
D(s_1,\dots ,s_k)_{\inv}$ and thus again by definition
$$r^M_U(\la(s_1,\dots ,s_k))=r^M_U(D(s_1,\dots ,s_k)).\qed$$
\enddemo

\proclaim{\nmb.{2.8}. Lemma }
The operator $\la$ commutes with Lie derivatives, i\.e\. for any
$X\in \Cal X(M)$ and any $s_i\in \Ga _c(E_iM)$
we have $\sum _{i=1}^k\la (s_1,\dots ,\Cal L_Xs_i,\dots ,s_k)=0$.
\endproclaim
\demo{Proof}
First we may clearly assume that $X$ has compact support. Thus we may
assume that the support of $X$ is arbitrarily small. So we may choose
$U_i$, $V_i$ and $f_i$ as in \nmb!{2.6} in such a way that the support of
$X$ is contained in $U_1$. Since $f\Cal L_Xs=\Cal L_X(fs)-(Xf)s$
and all $f_i$ are constant on the support of $X$, we see that
$\Cal L_X$ commutes with the multiplication by $f_i$ for each $i$
and thus we get
$$
\sum _{j=1}^k\la (s_1,\dots ,\Cal L_Xs_j,\dots ,s_k)
=\sum _{j=1}^k\sum_I
D(f_{i_1}s_1,\dots ,\Cal L_Xf_{i_j}s_j,\dots ,f_{i_k}s_k)_{\inv}.
$$
Now we exchange the two sums. For a fixed $I$ all sections
$f_{i_\ell}s_\ell$ and $\Cal L_Xf_{i_\ell}s_\ell$ vanish
simultaneously on some open subset. Then for the restrictions to
this subset we get
$$\multline
\sum _{j=1}^kD(f_{i_1}s_1,\dots ,\Cal L_Xf_{i_j}s_j,\dots ,
f_{i_k}s_k)_{\inv}=\\
=\sum _{j=1}^kD(f_{i_1}s_1,\dots ,\Cal L_Xf_{i_j}s_j,\dots ,f_{i_k}s_k)
=\Cal L_X(D(f_{i_1}s_1,\dots ,f_{i_k}s_k))=0,
\endmultline$$
since the restriction of $D(f_{i_1}s_1,\dots ,f_{i_k}s_k)$ to this
subset is absolutely invariant. Since absolutely invariant elements
are determined by their restrictions to any open subset the sum
vanishes globally and the result follows.
\qed\enddemo

\proclaim{\nmb.{2.9}. Theorem }
Let $E_1,\dots ,E_k$ be vector bundle functors defined on the category
$\Cal Mf_m^+$, $m\geq 2$, and let $\Cal F$ be an admissible natural
presheaf over a connected $m$--dimensional manifold $M$. Let
$D\:\Ga _c(E_1M)\x \dots \x \Ga _c(E_kM)\to \Cal F(M)$ be a $k$--linear
operator which commutes with Lie derivatives. Let
$\la\:\Ga _c(E_1M)\x \dots \x \Ga _c(E_kM)\to \Cal F(M)$ be the
$k$--linear operator constructed above, and put $\tilde D=D-\la$.
Then we have:
\roster
\item Both $\tilde D$ and $\la$ commute with Lie derivatives.
\item If $D$ is separately continuous then so are $\tilde D$
and $\la$.
\item $\la$ has values in the subspace of absolutely invariant
elements.
\item Suppose that $D$ satisfies the following condition:
\newline
Let $U\subset M$ be an open subset, $s_\ell\in \Ga_c(E_\ell M)$ for
$\ell=1,\dots ,k$. If there are $i$ and  $j$ between $1$ and $k$
such that $s_i$ vanishes on $U$ and the support of $s_j$ is contained
in $U$, then $D(s_1,\dots, s_k)$ restricts to zero on $U$.
\newline
Then $\la$ is weakly local and $\tilde D$ is local.
\endroster
\endproclaim
The only point that remains to be proved is (4). First we need a
lemma:

\proclaim{\nmb.{2.10}. Lemma }
Suppose that $D\:\Ga _c(E_1M)\x \dots \x \Ga _c(E_kM)\to \Cal F(M)$
satisfies the condition of \nmb!{2.9}.(4), and suppose that
$s_1,\dots ,s_k$ are sections such that for some $\ell\leq k$ and
some $i_1,\dots ,i_\ell$ between 1 and $k$ we have
$\supp (s_{i_1})\cap \dots \cap\supp (s_{i_\ell})=\emptyset$. Then
$D(s_1,\dots ,s_k)$ restricts to zero on
$M\setminus (\cup _{j=1}^\ell\supp s_{i_j})$.
\endproclaim
\demo{Proof}
We proceed by induction on $\ell$. For $\ell=1$ there is nothing to
prove. If $\ell=2$ we assume without loss of generality that
$s_1$ and $s_2$ have disjoint supports. Then $s_1$ vanishes on the
open set $M\setminus \supp(s_1)$ and $s_2$ has support contained in
this open set so by the condition \nmb!{2.9}.(4), $D(s_1,\dots ,s_k)$
restricts to zero on $M\setminus \supp(s_1)$ and thus also on
$M\setminus (\supp(s_1)\cup \supp(s_2))$.

So let us assume that $\ell> 2$ and that the result has been
proved for all integers smaller than $\ell$. Without loss of
generality we assume that
$\supp(s_1)\cap\dots\cap \supp(s_\ell)=\emptyset$ and that all
intersections of $\ell-1$ of these sets are nonempty. Then for
$i=1,\dots ,\ell$ we can construct open sets $U_i$ such that $U_i$ is
an open neighborhood of
$\supp(s_1)\cap\dots\cap\supp(s_{i-1})\cap
\supp(s_{i+1})\cap\dots\cap\supp(s_\ell)$, the sets $\bar U_i$ are
pairwise disjoint and $\bar U_i\cap\supp(s_i)=\emptyset$. Next for
any $i$ choose a smooth function $f_i\in C^\infty(M,[0,1])$ such that
$f_i$ has support contained in $U_i$ and is identically one on
$\supp(s_1)\cap\dots\cap\supp(s_{i-1})\cap
\supp(s_{i+1})\cap\dots\cap\supp(s_\ell)$, and put
$f_{\ell+1}=1-\sum_{i=1}^\ell f_i$. By multilinearity of $D$
we get $D(s_1,\dots ,s_k)=\sum_ID(f_{i_1}s_1,\dots ,f_{i_k}s_k)$,
where $I=(i_1,\dots ,i_k)$ runs over all $k$--tuples of integers
between $1$ and $\ell+1$. Now put
$W:=M\setminus (\cup _{i=1}^\ell\supp(s_i))$.
If there are at least two numbers between 1 and
$\ell$ in the multiindex $I$,
then the two corresponding functions have disjoint supports and
thus from above we see that the summand corresponding to $I$ vanishes
on an open set which certainly contains $W$. So
only summands corresponding to multiindices in which there is
only one fixed $j\leq \ell$ can contribute to the restriction to
$W$. But if $i_j=j$ then by
construction $f_{i_j}s_j=0$ and if $i_r=j$ for some $r\neq j$ then
the corresponding function has support disjoint to the support of
$s_j$ and so again the corresponding summand vanishes on an open set
which contains $W$. Thus we see that
$r^M_W(D(s_1,\dots ,s_k))=r^M_W(D(f_{\ell+1}s_1,\dots ,f_{\ell+1}s_k))$.
But by construction the support of $f_{\ell+1}$ is disjoint to
$\supp(s_1)\cap\dots\cap\supp(s_{\ell-1})$, so we see that
$\supp(f_{\ell+1}s_1)\cap\dots\cap\supp(f_{\ell+1}s_{\ell-1})=\emptyset$,
so by the induction hypothesis
$D(f_{\ell+1}s_1,\dots ,f_{\ell+1}s_k)$ vanishes on an open set which
contains $W$.
\qed\enddemo

\subheading{\nmb.{2.11} }
Next we prove that under the condition \nmb!{2.9}.(4) the operator
$\la$ is weakly local. So assume that $s_1,\dots ,s_k$ are sections
such that $\supp(s_1)\cap\dots\cap\supp(s_k)=\emptyset$. Choose any
triples $(U_i,V_i,f_i)$ as in \nmb!{2.6}. Then by \nmb!{2.10} for
any multiindex $I=(i_1,\dots ,i_k)$ the element
$D(f_{i_1}s_1,\dots ,f_{i_k}s_k)$ restricts to zero on the open
subset $M\setminus(\cup _j\supp(f_{i_j}s_j))$. By construction this set is
nonempty and moreover by definition the restrictions of
$D(f_{i_1}s_1,\dots ,f_{i_k}s_k)$ and
$D(f_{i_1}s_1,\dots ,f_{i_k}s_k)_{\inv}$ to this subset coincide.
Thus $D(f_{i_1}s_1,\dots ,f_{i_k}s_k)_{\inv}=0$ since absolutely
invariant elements are determined by their restrictions to arbitrary
open sets, and consequently $\la (s_1,\dots ,s_k)=0$.

\subheading{\nmb.{2.12} }
So it remains to show that under the condition \nmb!{2.9}.(4) the
operator $\tilde D$ is local. So assume that we have given sections
$s_1,\dots ,s_k$ and an open subset $W$ such that $s_1$ restricts to
zero on $W$. Then we have to show that
$r^M_W(D(s_1,\dots ,s_k))=r^M_W(\la(s_1,\dots ,s_k))$.
Let $U\subset W$ be open and such that $\bar U\subset W$.
Put $U_1:=U$, $V_1:=W$ and choose a function $f_1$
and triples $(U_i,V_i,f_i)$ for $i\geq2$ like in \nmb!{2.6}. Consider
a multiindex $I=(i_1,\dots ,i_k)$. If $i_1=1$ then $f_{i_1}s_1=0$ and
if $i_j=1$ for some other $j$ then the section
$f_{i_j}s_j$ has support contained in $W$, while $f_{i_1}s_1$
vanishes on $W$. So by the proof of \nmb!{2.10} and multilinearity we
see that $r^M_W(D(s_1,\dots ,s_k))=\sum _{I:i_j\neq 1}
r^M_W(D(f_{i_1}s_1,\dots ,f_{i_k}s_k))$.

Next $\la(s_1,\dots ,s_k)=
\sum _{I:i_j\neq 1}D(f_{i_1}s_1,\dots ,f_{i_k}s_k)_{\inv}$
by weak locality and multilinearity. For any
multiindex $I$ with all $i_\ell\neq 1$ the sections
$f_{i_1}s_1,\dots ,f_{i_k}s_k$ vanish simultaneously on the set
$U_1=U$. Consequently we see that
$r^M_U(\la (s_1,\dots ,s_k))=r^M_U(D(s_1,\dots ,s_k))$.

Now the set $W$ can be covered by countably many open subsets $W_i$
such that there are homeomorphisms $w_i:W_i\to \Bbb R^m$. For any
$(i,j)\in \Bbb N^2$ let $W_{ij}$ be the preimage under $w_i$ of the
open ball of radius $j$ around the origin. Then for all $i,j$ the
closure of $W_{ij}$ is contained in $W$. Next take a bijection
$\ka :\Bbb N\to \Bbb N^2$ and define
$U_\ell:=\cup _{i\leq \ell}W_{\ka (i)}$. Then clearly
$U_\ell\subset U_{\ell +1}$ for all $\ell$ and $\bar U_\ell\subset W$
and $W=\cup _{\ell\in \Bbb N}U_\ell$.

{}From above we see that $\tilde D(s_1,\dots ,s_k)$ restricts to zero
on each set $U_\ell$. Since the presheaf $\Cal F$ is admissible we
see that $\tilde D(s_1,\dots ,s_k)$ restricts to zero on $W$ by
condition (i) of \nmb!{2.1}.(3).

This finishes the proof of theorem \nmb!{2.9}.

\head\totoc\nmb0{3}. Weakly local natural functionals\endhead
In this section we give a complete description of weakly local
multilinear functionals which commute with Lie derivatives. This
classification can be found in \cite{Kirillov}, but
the proofs are only briefly sketched there. To have this paper self
contained we present full proofs here. Moreover we weaken some
continuity assumptions. The proofs are rather technical, so
readers who are not interested in details  should only read lemma
\nmb!{3.3} and skip the rest of this section.

\proclaim{\nmb.{3.1}. Lemma }
Let $\tau :C^\infty_c(\Bbb R^m,\Bbb R)\to \Bbb R$ be a not necessarily
continuous linear map such that for any
$f\in C^\infty_c(\Bbb R^m,\Bbb R)$ and any $i=1,\dots ,m$ we have
$\tau (\partial _if)=0$, where $\partial _if$ denotes the $i$--th
partial derivative of $f$. Then $\tau (f)=c\int_{\Bbb R^m}f(x)dx$ for
some real constant $c$. In particular $\tau$ is automatically
continuous.
\endproclaim
\demo{Proof}
A similar proof to that used for distributions in \cite{H\"ormander, 3.1.4
and 3.1.4'} applies. Let us first assume that $m=1$. Choose a function
$\ps \in C^\infty_c(\Bbb R,\Bbb R)$ such that $\int_{\Bbb R}\ps=1$.
Then for any $f\in C^\infty_c(\Bbb R,\Bbb R)$ we get
$$f(t)-\ps(t) \int_{\Bbb R}f=
\tfrac{d}{dt}\int_{-\infty}^t\biggl(f(x)-\ps (x)\int _{\Bbb R}f\biggr)dx,$$
and since the integral over the left hand side obviously vanishes,
the function which is differentiated on the right hand side has
compact support. Thus we get $0=\tau (f-\ps \int_{\Bbb R}f)$ and
hence $\tau (f)=\tau (\ps)\int_{\Bbb R}f$.

Next for $m>1$ take the function
$\ps$ from above, and for $g\in C^\infty_c(\Bbb R^{m-1},\Bbb R)$
define $g_m\in C^\infty_c(\Bbb R^m,\Bbb R)$ by
$g_m(x_1,\dots ,x_m):=g(x_1,\dots ,x_{m-1})\cdot \ps (x_m)$. Then
define $\tau _m:C^\infty_c(\Bbb R^{m-1},\Bbb R)\to \Bbb R$ by
$\tau _m(g)=\tau (g_m)$. For $f\in C^\infty_c(\Bbb R^m,\Bbb R)$ consider
$I_m(f)\in C^\infty_c(\Bbb R^{m-1},\Bbb R)$ defined by
$I_m(f)(x_1,\dots ,x_{m-1}):=\int _{\Bbb R}f(x_1,\dots ,x_m)dx_m$.
Now as above one shows that the function
$f(x_1,\dots ,x_m)-I_m(f)(x_1,\dots ,x_{m-1})\cdot \ps (x_m)$ is the
$m$--th partial derivative of a compactly supported smooth function
so we see that $\tau (f)=\tau ((I_m(f))_m)=\tau _m(I_m(f))$. Since by
construction $\tau_m$ is again linear an vanishes on partial
derivatives the result follows by induction.
\qed\enddemo

\proclaim{\nmb.{3.2}. Lemma}
Let $\la :(C^\infty_c(\Bbb R^m,\Bbb R))^k\to \Bbb R$ be a $k$--linear
functional of the form
$$\la (f_1,\dots ,f_k)=\sum _{\al^1,\dots ,\al^k}\int _{\Bbb R^m}
c_{\al^1\dots \al^k}(x)\partial^{\al^1}f_1(x)\dots \partial^{\al^k}f_k(x)dx,$$
where each $c_{\al^1\dots \al^k}$ is a continuous function on $\Bbb R^m$
and the sum is over finitely many $k$--tuples of multiindices
$\al^j=(\al^j_1,\dots ,\al^j_m)$ of nonnegative integers,
and $\partial^{\al^j}$ is the composition of partial derivatives
corresponding to $\al^j$. If $\la $ commutes with partial
derivatives, i\.e\.
$\sum _{i=1}^k\la (f_1,\dots ,\partial _\ell f_i,\dots ,f_k)=0$ for
all functions $f_j\in C^\infty_c(\Bbb R^m,\Bbb R)$ and any
$\ell =1,\dots ,m$, then there are constants $d_{\be^1\dots \be^{k-1}}$
such that
$$\la (f_1,\dots ,f_k)=\sum _{\be^1,\dots ,\be^{k-1}}\int _{\Bbb R^m}
d_{\be^1\dots \be^{k-1}}f_k(x)\partial^{\be^1}f_1(x)\dots
\partial^{\be^{k-1}}f_{k-1}(x)dx,$$
where again the sum is over a finite set of multiindices.
\endproclaim
\demo{Proof}
In order to bring a given expression for $\la$ to the required final
form we will rewrite its summands recursively. Therefore we need an
appropriate ordering on the sets of multi--multi indices. What we will
need is $(k-1)$--tuples $(\al^1,\dots ,\al^{k-1})$ where each
$\al ^j=(\al^j_1,\dots ,\al^j_m)$ is an $m$--tuple of nonnegative
integers. First we order these $m$--tuples  by their total degree
$|\al ^j|=\sum \al^j_i$ and within each total degree
lexicographically. Then we order the $(k-1)$--tuples by their total
degree $|(\al^1,\dots ,\al^{k-1})|=\sum |\al ^j|$ and within each total
degree lexicographically.

So let assume that for some fixed $(k-1)$--tuple
$(\be^1,\dots ,\be^{k-1})$ and fixed $m$--tuple $\be$ we already have
an expression for $\la$ of the following form:
$$\align
\la (f_1,\dots ,f_k)&= \kern-25pt\sum_{(\al^1,\dots ,\al^{k-1})<(\be^1,\dots
,\be^{k-1})}
\int_{\Bbb R^m}d_{\al^1\dots \al^{k-1}}f_k(x)\partial ^{\al^1}f_1(x)\dots
\partial ^{\al^{k-1}}f_{k-1}(x)dx+\\
&+\sum_{\be^k\leq \be}\int_{\Bbb R^m}c_{\be^1\dots \be^k}(x)
\partial ^{\be^1}f_1(x)\dots \partial ^{\be^k}f_k(x)dx+\\
&+\sum_{(\al^1,\dots ,\al^{k-1})>(\be^1,\dots ,\be^{k-1})}\sum_{\al^k}
\int_{\Bbb R^m}c_{\al^1\dots \al^k}(x)\partial ^{\al^1}f_1(x)\dots
\partial ^{\al^k}f_k(x)dx,\tag{1}
\endalign$$
where the $d$'s are constants, the $c$'s are continuous functions and
$c_{\be^1\dots \be^{k-1}\be}$ is nonzero.

Now fix a relatively compact open subset $U\subset \Bbb R^m$. Then choose
a compactly supported smooth function $f$ on $\Bbb R^m$ such that
$f$ is identically one on $U$, and consider the distribution
$f_k\mapsto \la (x^{\be ^1}f,\dots ,x^{\be^{k-1}}f,f_k)$ on
$C^\infty_c(U,\Bbb R)$. Then this is given by
$$\align
f_k\mapsto\ &\sum\Sb(\al^1,\dots ,\al^{k-1})<(\be^1,\dots ,\be^{k-1})\\
\al^i\le\be^i\endSb
\int_{\Bbb R^m}d_{\al^1\dots \al^{k-1}}\tfrac{\be^1!}{(\be^1-\al^1)!}
\dots \tfrac{\be^{k-1}!}{(\be^{k-1}-\al^{k-1})!}f_k(x)\cdot\\
&\hbox to2in{\hfil}\cdot x^{\be^1-\al^1}\dots x^{\be^{k-1}-\al^{k-1}}dx+\\
&+\sum_{\be^k\leq \be}\int_{\Bbb R^m}\be^1!\dots \be^k!
c_{\be^1\dots \be^k}(x)\partial ^{\be^k}f_k(x)dx\tag{2}
\endalign$$
where the first sum collects all summands with
$(\al^1,\dots ,\al^{k-1})<(\be^1,\dots ,\be^{k-1})$, since the
remaining terms vanish.
Using the commutation of $\la$ with partial derivatives and (1) we
get:
$$\multline
\la (x^{\be ^1}f,\dots ,x^{\be^{k-1}}f,\partial _if_k)=
-\sum_{j=1}^{k-1}\la (x^{\be ^1}f,\dots ,\partial _i(x^{\be^j}f),\dots
,x^{\be^{k-1}}f,f_k)=\\
=-\sum_j\sum\Sb(\al^1,\dots ,\al^{k-1})<(\be^1,\dots ,\be^{k-1})\\
\al^i\le\be^i\endSb
\int_{\Bbb R^m}d_{\al^1\dots \al^{k-1}}\tfrac{\be^1!}{(\be^1-\al^1)!}
\dots \tfrac{\be^{k-1}!}{(\be^{k-1}-\al^{k-1})!}f_k(x)\cdot\\
\cdot x^{\be^1-\al^1}\dots
\partial_i(x^{\be^j-\al^j})\dots x^{\be^{k-1}-\al^{k-1}}dx+0.
\endmultline\tag3$$
On the other hand using (2) we compute:
$$\align
\la (x^{\be ^1}f,&\dots ,x^{\be^{k-1}}f,\partial _if_k)=\\
&\sum\Sb(\al^1,\dots ,\al^{k-1})<(\be^1,\dots ,\be^{k-1})\\\al^i\le\be^i
\endSb
\int_{\Bbb R^m}d_{\al^1\dots \al^{k-1}}\tfrac{\be^1!}{(\be^1-\al^1)!}
\dots \tfrac{\be^{k-1}!}{(\be^{k-1}-\al^{k-1})!}\partial_if_k(x)\cdot\\
&\hbox to2in{\hfil}\cdot x^{\be^1-\al^1}\dots x^{\be^{k-1}-\al^{k-1}}dx+\\
&+\sum_{\be^k\leq \be}\int_{\Bbb R^m}\be^1!\dots \be^k!
c_{\be^1\dots \be^k}(x)\partial_i\partial ^{\be^k}f_k(x)dx.\tag4
\endalign$$
Now all terms in the first sum of (4) can be integrated by parts and
this gives exactly the sum in (3). Thus the remaining sum must be
zero. This means that for the restrictions of the $c$'s to $U$ we get
the equation
$$0=\partial_i\bigl(\sum_{\be^k\leq \be}(-1)^{|\be^k|+1}\be^1!\dots \be^k!
\partial ^{\be^k}c_{\be^1\dots \be^k}\bigr),\tag5$$
where the derivatives are interpreted in the distributional sense.
Since this equation holds over each relatively compact open subset of
$\Bbb R^m$ it holds on $\Bbb R^m$.

Assume first that $\be =(0,\dots ,0)$. Then (5) implies that
$c_{\be^1\dots \be^{k-1}\be}$ is constant on $\Bbb R^m$ and thus we
can rewrite the expression (1) for $\la$ with $(\be^1,\dots ,\be^{k-1})$
replaced by some greater $(\ga^1,\dots ,\ga^{k-1})$. Moreover in
this procedure we do not increase the maximal total degree occurring
in the expression.

Next assume that $\be >(0,\dots ,0)$. Then we see from (5) that there
is a constant $C$ such that
$$(-1)^{|\be|+1}\be !\partial ^{\be}c_{\be^1\dots \be^{k-1}\be}=
C+\sum_{\be^k<\be}(-1)^{|\be^k|}\be^k!\partial^{\be^k}c_{\be^1\dots \be^k}.$$
Thus for arbitrary $f_1,\dots ,f_k\in C^\infty_c(\Bbb R^m,\Bbb R)$ we
have
$$\multline
-\int_{\Bbb R^m}\be !c_{\be^1\dots \be^{k-1}\be}(x)\partial ^{\be}\left(
f_k\partial^{\be^1}f_1\dots \partial^{\be^{k-1}}f_{k-1}
\right)(x)dx=\\
=C\int_{\Bbb R^m}f_k(x)\partial^{\be^1}f_1(x)\dots
\partial^{\be^{k-1}}f_{k-1}(x)dx+\\
+\sum_{\be^k<\be}\be^k!c_{\be^1\dots \be^k}(x)\partial ^{\be^k}\left(
f_k\partial^{\be^1}f_1\dots \partial^{\be^{k-1}}f_{k-1}
\right)(x)dx.
\endmultline$$
Expanding the partial derivatives of products in this identity we see
that we can express
$$\int_{\Bbb R^m}c_{\be^1\dots \be^{k-1}\be}(x)
\partial^{\be^1}f_1(x)\dots \partial^{\be^{k-1}}f_{k-1}(x)
\partial ^{\be}f_k(x)dx$$
as a sum of terms of the form
$\int_{\Bbb R^m}\tilde c_{\ga^1\dots \ga^k}(x)
\partial^{\ga^1}f_1(x)\dots \partial ^{\ga^k}f_k(x)dx$
in which either $(\ga^1,\dots ,\ga^{k-1})>(\be^1,\dots ,\be^{k-1})$
or $(\ga^1,\dots ,\ga^{k-1})=(\be^1,\dots ,\be^{k-1})$ and
$\ga^k<\be$. Moreover we always have
$|\ga^1|+\dots +|\ga^k|\leq |\be^1|+\dots +|\be^{k-1}|+|\be|$.

Thus finitely many applications of this procedure lead to an
expression for $\la$ of the required form.
\qed\enddemo

\proclaim{\nmb.{3.3}. Lemma }
Let $E_1,\dots ,E_k$ be natural vector bundles defined on
$\Cal Mf_m^+$, and let $\la\:\Ga _c(E_1M)\x\dots\x\Ga_c
(E_kM)\to \Bbb R$ be a weakly local separately continuous
$k$--linear map, where $M$ is a connected $m$--dimensional oriented
manifold. Then there is a $(k-1)$--linear local natural operator
$D^{\la}\: E_1\x\dots\x E_{k-1}\to E_k^*\otimes \La^mT^*$, where
$E_k^*$ is the bundle functor which assigns to each manifold $M$ the dual
bundle of $E_kM$, such that
$$\la (s_1,\dots ,s_k)=\int_M\langle D^{\la}(s_1,\dots ,s_{k-1}),s_k\rangle .$$
Here $\langle \quad,\quad\rangle $ denotes the canonical pairing
$$\Ga_c(E_k^*M\otimes \La^mT^*M)\x\Ga_c(E_kM)\to \Ga_c(\La^mT^*M).$$
In particular $\la$ is automatically jointly continuous.

If $k=1$, i\.e\. $\la$ is linear, this means that there is an
absolutely invariant section $\si$ of $E_1^*M\otimes\La^mT^*M$ such
that $\la (s)=\int _M\langle \si, s\rangle$. Moreover in this case
the continuity follows from the other assumptions.
\endproclaim
\demo{Proof}
Let us first consider the case $M=\Bbb R^m$. Up to a natural
isomorphism, for each natural bundle $E_j$ we have
$E_j\Bbb R^m=\Bbb R^m\x V_j$,
with the action of the translations
$t_x\: y\mapsto x+y$ given by $t_x \x \text{id}_{V_j}$. Thus the Lie
derivatives with respect to the constant fields are just the partial
derivatives of the coordinate functions. (This identification is
always achieved by $E_j\Bbb R^m \ni v\mapsto (p_j(v),
E(t_{-p_j(v)})(v))\in \Bbb R^m\x V_j$.)

Now in any $V_j$ choose a fixed basis. Then for any multiindex
$(i_1,\dots ,i_k)$ with $1\leq i_j\leq\text{dim}(V_j)$ consider the
$k$--linear functional
$\la_{i_1\dots i_k}:C^\infty(\Bbb R^m,\Bbb R)^k\to \Bbb R$ which gives
the coordinate expression of $\la$. These functionals are clearly
separately continuous,
weakly local and they commute with partial derivatives if $\la$ has
these properties.

Let us first treat the case $k=1$, i\.e\. $\la$ is linear. Then
without continuity assumptions we conclude from lemma \nmb!{3.1} that
for each $i_1$ there is a constant $c_{i_1}$ such that
$\la_{i_1}(f)=c_{i_1}\int_{\Bbb R^m}f(x)dx$. Hence for any section
$f=(f^1,\dots ,f^{n_1})\in C^\infty_c(\Bbb R^m,V_1)$ we have
$\la(f)=\int_{\Bbb R^m}\sum_ic_if^i(x)dx$. Now we interpret
$(c_1,\dots ,c_{n_1})dx$ as a section $\si$ of
$\Bbb R^m\x V_1^*\otimes \La^m\Bbb R^{m*}=E_1^*\Bbb R^m\otimes
\La^mT^*\Bbb R^m$. Hence we have
$\la (f)=\int_{\Bbb R^m}\langle f,\si\rangle$. Since any Lie
derivative of an $m$--form on $\Bbb R^m$ is exact we get:
$$\multline
0=\int_{\Bbb R^m}\Cal L_X(\langle \si, f\rangle )=
\int_{\Bbb R^m}\langle \si,\Cal L_Xf\rangle +
\int_{\Bbb R^m}\langle \Cal L_X\si, f\rangle =\\
=\la (\Cal L_Xf)+\int_{\Bbb R^m}\langle \Cal L_X\si, f\rangle
=\int_{\Bbb R^m}\langle \Cal L_X\si, f\rangle .
\endmultline$$
Since this holds for all $f\in \Ga _c(EM)$ and
$X\in \Cal X(\Bbb R^m)$ we have
$\Cal L_X\si=0$ and thus $\si$ is absolutely invariant.

So let us turn back to the general multilinear case. Applying the
Schwartz kernel theorem (c\.f\. \cite{H\"ormander, 5.2.1}) to
$\la_{i_1\dots i_k}$ we see that this functional is given by a
distribution on $(\Bbb R^m)^k$, which has support in the diagonal by
the weak locality. Thus it is given (at least locally) as
$$\la_{i_1\dots i_k}(f_1,\dots ,f_k)=\sum_{\al_1,\dots ,\al_k}
\int_{\Bbb R^m}c_{\al_1\dots \al_k}(x)\partial^{\al_1}f_1(x)\dots
\partial^{\al_k}f_k(x)dx,$$
where the $c$'s are continuous functions and the sum is over some
finite number of multiindices $\al_j$ (c\.f\. \cite{H\"ormander, 5.2.3}).

By lemma \nmb!{3.2} we can express $\la_{i_1\dots i_k}$ locally as
$$\sum _{\be^1,\dots ,\be^{k-1}}\int _{\Bbb R^m}
d_{\be^1\dots \be^{k-1}}^{i_1\dots i_k}f_k(x)\partial^{\be^1}f_1(x)\dots
\partial^{\be^{k-1}}f_{k-1}(x)dx.$$
Fix an open subset $U$ over which such a representation is valid.
Then for smooth sections $f_j\in C^\infty_c(U,V_j)$ we have
$$\la(f_1,\dots ,f_k)=\sum_{i_1,\dots ,i_k}\sum _{\be^1,\dots ,\be^{k-1}}
\int _{\Bbb R^m}d_{\be^1\dots \be^{k-1}}^{i_1\dots i_k}f_k^{i_k}(x)
\partial^{\be^1}f_1^{i_1}(x)\dots \partial^{\be^{k-1}}f_{k-1}^{i_{k-1}}(x)dx.$$
We can now define an operator
$D_U:C^\infty_c(U,V_1)\x \dots \x C^\infty_c(U,V_{k-1})\to
C^\infty_c(U,V_k^*\otimes \La^m\Bbb R^{m*})$  by
$$\langle D_U(f_1,\dots ,f_{k-1}),f_k\rangle:=
\sum \Sb \be^1,\dots ,\be^{k-1}\\{i_1,\dots ,i_k}\endSb
d_{\be^1\dots \be^{k-1}}^{i_1\dots i_k}f_k^{i_k}(x)
\partial^{\be^1}f_1^{i_1}(x)\dots
\partial^{\be^{k-1}}f_{k-1}^{i_{k-1}}(x)dx.$$
Clearly this is a well defined local operator and over $U$,
$$\la(f_1,\dots ,f_k)=
\int_{\Bbb R^m}\langle D_U(f_1,\dots ,f_{k-1}),f_k\rangle,$$
and $D_U$ is uniquely determined by this formula.
As in the linear case we compute now:
$$\align
0&=\int_{\Bbb R^m}\Cal L_X\langle D_U(f_1,\dots ,f_{k-1}),f_k\rangle=\\
&=\int_{\Bbb R^m}\langle\Cal L_X D_U(f_1,\dots ,f_{k-1}),f_k\rangle+
\int_{\Bbb R^m}\langle D_U(f_1,\dots ,f_{k-1}),\Cal L_Xf_k\rangle=\\
&=\int_{\Bbb R^m}\langle\Cal L_X D_U(f_1,\dots ,f_{k-1}),f_k\rangle+
\la(f_1,\dots ,f_{k-1},\Cal L_Xf_k)=\\
&=\int_{\Bbb R^m}\langle\Cal L_X D_U(f_1,\dots ,f_{k-1}),f_k\rangle-
\sum_{i=1}^{k-1}\la(f_1,\dots ,\Cal L_Xf_i,\dots ,f_k)=\\
&=\int_{\Bbb R^m}\langle\Cal L_X D_U(f_1,\dots ,f_{k-1}),f_k\rangle-
\sum_{i=1}^{k-1}\langle D_U(f_1,\dots ,\Cal L_Xf_i,\dots ,f_{k-1}),f_k\rangle .
\endalign$$
Thus $D_U$ commutes with Lie derivatives, so it extends uniquely to a
natural operator by \cite{\v Cap--Slov\'ak, 3.2}, i\.e\. there is a unique
natural
operator $D$ on $\Cal Mf_m^+$ such that its value on $U$ is $D_U$.

Any point of $\Bbb R^m$ has a neighborhood on which the construction
can be carried out as above. Moreover from the uniqueness of the
extension we see that for intersecting open sets the natural
operators obtained coincide, and thus $\la$ is given on $\Bbb R^m$ by
the formula
$$\la(f_1,\dots ,f_k)=
\int_{\Bbb R^m}\langle D_{\Bbb R^m}(f_1,\dots ,f_{k-1}),f_k\rangle.$$

Finally take an arbitrary manifold $M$. Choose an oriented atlas
$(U_i,u_i)$ of $M$ such that any $u_i$ is a diffeomorphism onto
$\Bbb R^m$. Then we have the induced isomorphisms
$E_j(u_i):E_jM\supset E_jU_i\to E_j\Bbb R^m$. Any section of $E_jM$
with support in $U_i$ can be written as $u_i^*f_j$ for a compactly
supported section $f_j$. Now by the first part of this proof there is
a unique natural operator $D$ such that
$$\align
\la(u_i^*f_1,\dots ,u_i^*f_k)&=\int_{\Bbb R^m}\langle
D_{\Bbb R^m}(f_1,\dots ,f_{k-1}),f_k\rangle=\\
&=\int_{\Bbb R^m}\langle (u_i^{-1})^*D_M(u_i^*f_1,\dots ,u_i^*f_{k-1}),f_k
\rangle=\\
&=\int_M\langle D_M(u_i^*f_1,\dots ,u_i^*f_{k-1}),u_i^*f_k\rangle .
\endalign$$
As before the uniqueness implies that the operators obtained from
intersecting charts coincide, and since $M$ is connected and any
compactly supported section can be written as a finite sum of
sections having support in some chart $U_i$ the proof is finished.
\qed\enddemo

\head\totoc\nmb0{4}. The linear case\endhead
Now we give a complete description of linear
operators between sections of natural vector bundles over oriented
manifolds which commute
with Lie derivatives. In the next section this result will be used as the
starting point for an induction procedure.

\proclaim{\nmb.{4.1}. Proposition }
Let $F$ be a natural vector bundle defined on the category
$\Cal Mf_m^+$ of oriented $m$--dimensional manifolds and orientation
preserving local diffeomorphisms
and let $M$ be such a manifold which is connected.
Then the sheaf of smooth sections of
$FM$ is an admissible natural presheaf over $M$.
\endproclaim
\demo{Proof}
By the general theory mentioned in \nmb!{1.2}, there is the
canonically defined continuous action of the sheaf of Lie algebras of vector
fields on the base manifold on the sections of the natural vector bundle.

Since the sections of a vector bundle form a sheaf, condition (i) of
\nmb!{2.1}.(3)  is trivially satisfied. On the other hand the
existence and uniqueness condition in
(ii) of \nmb!{2.1}.(3) follows immediately from \nmb!{1.3}. Finally
continuity of the extension operators is trivial since the absolutely
invariant sections always form a finite dimensional space.
\qed\enddemo

\subheading{\nmb.{4.2} }
Let $E$ and $F$ be natural vector bundles over $\Cal Mf_m^+$, let $M$
be a connected oriented manifold of dimension $m\geq 2$ and let
$D:\Ga _c(EM)\to \Ga (FM)$ be a linear operator which
commutes with Lie derivatives. Then the condition of \nmb!{2.9}.(4)
is automatically satisfied since in the linear case it just means
that $D(0)=0$.
Thus we may apply theorem \nmb!{2.9} to split $D$ into a local linear
operator $\tilde D$ and a linear map $\la$ with values in the space
of absolutely invariant sections which vanishes on Lie derivatives.
Theorem 3.2 of \cite{\v Cap--Slov\'ak} shows that $\tilde D$ is a natural
operator, which in particular implies that it is continuous.

To discuss the linear map $\la$ let $v_1,\dots ,v_n$ be a basis of
the space of absolutely invariant sections of $FM$ and let $\{
v_i^*\}$ be the dual basis. Then for any $i=1,\dots ,n$ the map
$\la ^i:=v_i^*\o \la :\Ga _c(EM)\to \Bbb R$ is a linear
functional which vanishes on Lie derivatives. Using lemma \nmb!{3.3}
we get:

\proclaim{\nmb.{4.3}. Theorem }
Let $M$ be a connected oriented smooth manifold of dimension
$m\geq 2$, $E$ and $F$ natural vector bundles defined over
$\Cal Mf_m^+$, $v_1,\dots ,v_n$ a basis for the space of absolutely
invariant sections of $FM$ and let $D:\Ga _c(EM)\to \Ga (FM)$ be a
linear operator which commutes with Lie derivatives. Then
there is a (local) natural operator $\tilde D:\Ga (EM)\to \Ga (FM)$
and there are absolutely invariant sections
$\si_1,\dots ,\si_n\in \Ga (E^*M\otimes \La ^mT^*M)$ such that
$$
D(s)=\tilde D(s)+\sum _{i=1}^n(\int_M\langle \si_i, s\rangle )v_i.
$$
So $D$ is an almost natural operator and thus in particular
automatically continuous.
\endproclaim

\proclaim{\nmb.{4.4}. Corollary }
In the setting of theorem \nmb!{4.3} assume that either $FM$ or
$E^*M\otimes \La ^mT^*M$ has no absolutely invariant sections.
Then any linear operator $\Ga _c(EM)\to \Ga (FM)$ or
$\Ga (EM)\to \Ga (FM)$ which commutes with Lie derivatives is a
local natural operator.
\endproclaim
Note that the conditions in the corollary are satisfied in many
concrete situations. In subbundles of tensor bundles for example
absolutely invariant sections can only exist if the numbers of
covariant and contravariant indices are equal.

\proclaim{\nmb.{4.5}. Corollary }
In the setting of theorem \nmb!{4.3} assume that $M$ is not compact.
Then any linear operator $\Ga _c(EM)\to \Ga _c(FM)$ as
well as any continuous linear operator $\Ga(EM)\to \Ga(FM)$ which
commutes with Lie derivatives is a (local) natural operator.
\endproclaim
\demo{Proof}
In the case of compact supports the result follows immediately from
the fact that on a noncompact manifold no nonzero absolutely invariant
section can have compact support by \nmb!{4.1}.
In the second case we may restrict the operator to an operator
$\Ga _c(EM)\to \Ga (FM)$. There it splits by theorem \nmb!{4.3} as
$\tilde D+\sum_i\la_iv_i$, where the $v_i$ are a basis for the space
of absolutely invariant sections of $FM$. Now the natural operator
$\tilde D$ extends continuously to $\Ga (EM)$, and we consider the
difference $\la$ between this extension and the original operator. On
compactly supported sections this difference has the form
$\sum_i\la_iv_i$, where $\la_i(s)=\int_M\langle \si_i, s\rangle $, for an
absolutely
invariant section $\si_i\in \Ga (E^*M\otimes \La ^mT^*M)$. Now assume
that some $\si _i$, say $\si _1$ is nonzero. Then using a chart
construction one easily concludes that for any open subset
$U\subset M$ there is a section $s\in \Ga_c(EM)$ with support
contained in $U$ such that $\la_1(s)=1$.

Now choose compact subsets $K_i\subset M$ such that each $K_i$ is
contained in the interior of $K_{i+1}$ and such that
$M=\cup_{i\in \Bbb N}K_i$. Then for any $i\in \Bbb N$ choose an open
subset $U_i\subset K_i\setminus K_{i-1}$ and a smooth section
$s_i\in \Ga _c(EM)$ with support contained in $U_i$ such that
$\la _1(s_i)=1$. Then any compact subset of $M$ is contained in some
$K_n$ and thus intersects only finitely many $U_i$. Thus the sum
$\sum_{i\in \Bbb N}s_i$ is actually finite over each compact and thus
converges in $\Ga (EM)$. By continuity of $\la$ we see that this
element is mapped by $\la$ to the limit of $\sum _{i=1}^k\la (s_i)$.
But this sum is divergent by construction.

Thus we see that $\la$ vanishes on all compactly supported sections
and since these are dense in all sections the result follows.
\qed\enddemo

\proclaim{\nmb.{4.6}. Corollary}
Let $E$ and $F$  be vector bundle functors on $\Cal Mf_m^+$.
The natural presheaf $\Cal F(U):=L(\Ga _c(EU);\Ga (FU))$
of continuous linear operators is admissible.
\endproclaim
\demo{Proof}
By definition, an element $D\in\Cal F(U):=L(\Ga _c(EU);\Ga (FU))$
is absolutely invariant if and only if $D$ commutes with Lie
derivatives. But then $D$ is of the form
$D(s)=\tilde D(s) + \sum_i(\int_M\langle \si_i,s\rangle )v_i$ with $\tilde D$
natural and $\si_i$ absolutely invariant. Thus both summands are
completely determined by any restriction to an open submanifold.
Moreover using \nmb!{1.6} and \nmb!{4.1} we see that the space of
absolutely invariant elements is finite dimensional and thus condition
(ii) of \nmb!{2.1}.(3) is satisfied.

Let $U\subset M$ be open and let $U_i$, $i\in \Bbb N$, be a system of
open subsets such that $U_i\subset U_{i+1}$ for all $i$ and such that
$U=\cup _{i\in \Bbb N}U_i$.
Suppose that $D\in L(\Ga _c(E_U);\Ga (FU))$ restricts to zero on each
$U_i$. Take a section $s\in \Ga _c(EU)$. Then the support of $s$ is
contained in some $U_i$ and thus in all $U_j$ where $j>i$. Hence by
assumption $D(s)$ vanishes on each $U_j$ and thus on the union of all
$U_j$ and so condition (i) of \nmb!{2.1}.(3) is satisfied as well.
\qed\enddemo

\head \nmb0{5}. Proof of Theorem \nmb!{1.9}\endhead

The main idea is to view $k$--linear operators between sections of
vector bundles as $\ell$--linear operators with
values in the natural presheaf of $(k-\ell)$--linear operators,
as defined in example \nmb!{2.3}.(2), to apply theorem \nmb!{2.9} and
to use induction.

\subheading{\nmb.{5.1}. Definition}
Let $E_1,\dots ,E_k,E$ be natural vector bundles defined on
$\Cal Mf_m^+$, and let $D:\Ga _c(E_1M)\x\dots\x\Ga_c
(E_kM)\to \Ga(EM)$ be a $k$--linear separately continuous operator,
where $M$ is an $m$--dimensional oriented manifold. For an
$\ell$--tuple $(i_1,\dots ,i_\ell)$ of integers with
$1\leq i_1<\dots <i_\ell\leq k$ we define the {\sl associated
$\ell$--linear operator of type $(i_1,\dots ,i_\ell)$ to $D$\/} to be
the operator
$$D^{i_1,\dots ,i_\ell}:\Ga_c(E_{i_1}M)\x\dots\x\Ga_c(E_{i_\ell}M)\to
L(\Ga_c(E_{i_{\ell +1}}M)\x\dots\x\Ga_c(E_{i_k}M);\Ga (EM))$$
given by
$D^{i_1,\dots ,i_\ell}(s_{i_1},\dots ,s_{i_\ell})
(s_{i_{\ell +1}},\dots ,s_{i_k}):=D(s_1,\dots ,s_k)$.
Here $(i_{\ell +1},\dots ,i_k)$ denotes the
ordered sequence of integers between 1 and $k$ which do not occur in
$(i_1,\dots ,i_\ell)$.

\subheading{\nmb.{5.2} }
In the notation of \nmb!{5.1} assume that $D$ commutes with Lie
derivatives. Then by definition of the Lie derivatives on the natural
presheaf of $(k-\ell)$--linear operators any associated $\ell$--linear
operator commutes with Lie derivatives, too.

Next the condition \nmb!{2.9}.(4) can be interpreted nicely in this
situation.  Consider the operator $D^{i_1,\dots ,i_\ell}$ from above.
Then the condition of \nmb!{2.9}(4) for fixed $j$ means that if any
of the other sections $s_{i_1},\dots ,s_{i_\ell}$ vanishes on an open
subset $U\subset M$ which contains the support of $s_{i_j}$ then
$D^{i_1,\dots ,i_\ell}(s_{i_1},\dots ,s_{i_\ell})$ restricts to zero
on $U$. This in turn means that if all $s_{i_{\ell+1}},\dots ,s_{i_k}$
have support contained in $U$, then $D(s_1,\dots ,s_k)$ restricts
to zero on $U$. But
this is exactly equivalent to locality of the operator
$D^{i_1,\dots ,i_{j-1},i_{j+1},\dots ,i_\ell}$.

Thus the condition in \nmb!{2.9}.(4) means for $\ell$--linear
operators with values in the presheaf of $(k-\ell)$--linear operators
of the type we consider exactly that all associated
$(\ell-1)$--linear operators with values in the presheaves of
$(k-\ell+1)$--linear operators are local.

Before we can prove our main theorem we have to study the locality
properties of almost natural operators.

\proclaim{\nmb.{5.3}. Lemma}
Let $D$ be an elementary almost natural $k$--linear operator of the type
$(I^1,\dots ,I^r,J)$,
and  let $I$ be any subset of $\{ 1,\dots ,k\}$. Then the associated
operator $D^I$ is local if and only if no $I^j$ is a subset of $I$.
\endproclaim
\demo{Proof}
Let us assume that no $I^j$ is a subset of $I$, take some $i\in I$
and consider a section $s_i$ which vanishes on some open subset $U$
of $M$. If $i\in J$ then clearly $D(s_1,\dots ,s_k)$ vanishes on $U$
for arbitrary sections $s_\ell$, $\ell\neq i$. On the other hand if
$i$ is in some $I^j$ then by assumption there is some $\ell\in I^j$
which is not in $I$. Then for any section $s_\ell$ with support
contained in $U$ we have $\la ^j(\dots ,s_i,\dots ,s_\ell,\dots)=0$
by weak locality of $\la ^j$.

Conversely assume that $I^j\subset I$. Choose a small open subset $U$
of $M$. Then we can find sections $s_n$ for $n\in I^j$ which have
support in $M\setminus U$ such that $\la^j$ is nonzero on these
sections. On the other hand for all $\ell \neq j$ we can find
sections corresponding to $I^\ell$ with support contained in $U$ on
which $\la ^\ell$ is nonzero. Moreover we can also find sections
corresponding to $J$ with support in $U$ on which the local natural
part does not vanish, and hence we get a nonzero result over $U$
which contradicts locality of $D^I$.
\qed\enddemo

\subheading{\nmb.{5.4}}
Next let us consider a general almost natural operator $D$. By
definition we can write $D$ as a sum of nonzero elementary almost
natural operators.

\proclaim{Lemma}
Let $D$ be an almost natural operator
and let $I$ be a subset of $\{ 1,\dots ,k\}$. Then the associated
operator $D^I$ is local if and only if the associated operator of
type $I$ to any of the elementary summands is local.

Moreover, any set of elementary  almost natural operators of pairwise
different type is linearly independent.
\endproclaim
\demo{Proof}
The non-trivial part of the proof is to show that if $D^I$ is local
then each elementary summand has the same locality property. By
\nmb!{5.3} this means that the elementary summands  must not involve a
set $I^j\subset I$ in their types.
We may assume without loss of generality
that if there are several elementary almost natural
operators of the same type in this sum, then their local natural parts
(corresponding to the last set in the type) are linearly independent.
{}From now on we will assume that all representations of almost natural
operators as sums of elementary ones are of this form.

Suppose that $D$ and $I$ are given and take some $i_0\in I$ such that
there is an elementary summand which has in its type a set
$I^j\subset I$ which contains $i_0$, and consider only those summands
which have this property. Next consider only those operators for
which the corresponding set $I^j$ is of minimal cardinality. Next
take the minimal cardinality of the other sets in the types of the
remaining operators and consider only those which have a set of this
cardinality in their type. Next take again the minimal cardinality of
the remaining sets and the corresponding operators and so on up to a
point where there are only operators left which have all the same
cardinalities of sets in their types. Now choose one of the remaining
operators and renumber the bundles in such a way that $i_0$ becomes 1
and the operator has type $I^1=\{ 1,2,\dots ,n_1\} ,\dots ,I^r=\{
n_1+\dots +n_{r-1}+1,\dots ,n_1+\dots +n_r\} ,J=\{ n_r+1,\dots ,k\}$,
with $n_2\le n_3\le\dots\le n_r$.
Assume that it is of the form $\la^1\dots \la ^r\tilde D$.
Now choose an open subset $U$ of $M$ and sections
$s_1,\dots ,s_{n_1}$ which all vanish on $U$ such that $\la ^1$ is
nonzero on these sections. Next choose $r$ disjoint open subsets
$U_2,\dots ,U_{r+1}$ contained in $U$, for any $i\le r$ choose sections with
support in $U_i$ on which $\la ^i$ is nonzero and choose sections
with support contained in $U_{r+1}$ on which $\tilde D$ is nonzero.

Since all $\la^i$ are weakly local (cf\. \nmb!{1.7})
we see that any operator in the linear
combination which does not vanish identically on these sections must
belong to those which remained at the end of our choices.
So let us consider such operators, and let us denote their types with tildes.
By definition of the type we must have $1\in \tilde I^1$ or
$1\in \tilde J$. By minimality of the cardinality of $I^1$ we must
then have $\tilde I^1=I^1$ or $\tilde J=I^1$ to get something
nonzero. If $\tilde J=I^1$ we must have $n_1+1\in \tilde I^1$ so
again by minimality we must have $\tilde I^1=I^2$, then
$\tilde I^2=I^3$ and so on. On the other hand if $\tilde I^1=I^1$
then $n_1+1\in \tilde I_2$ or $n_1+1\in \tilde J$ and again the only
possibility to get a different type is $\tilde J=I^2$,
$\tilde I^2=I^3$ and so on.
Thus we see that the only operators which can produce something
nonzero are those whose types satisfy $\tilde I^j=I^j$ for $j\leq\ell$,
$\tilde I^j=I^{j+1}$ for $\ell<j\le r$, $\tilde I^r=J$  and
$\tilde J=I^{\ell+1}$, or those with $\tilde J=\emptyset$ for which
we must have $\tilde I^j=I^j$ for $j\le r$ and $\tilde I^{r+1}=J$.

If $D^I$ is local, these operators must add up to zero on $U$ for
each choice of sections as above.
The operators with $\tilde J\neq \emptyset$ produce section with
support in $U_\ell$ if $\tilde J=I^\ell$, while the ones with
$\tilde J=\emptyset$ produce absolutely invariant sections which in
particular have no zeros. Thus for any remaining type the operators
of this type have to add up to zero on some $U_\ell$ or on $U$.
Since we may assume that
for each fixed type the occurring local natural parts are linearly
independent this implies that the individual elementary summands must
already be zero.
Thus, altogether, if $D^I$ is local, then there is no elementary
summand with $I_j\subset I$ in its type.

This concludes the proof of the first assertion of the
lemma. The linear independence of any set
of elementary almost natural operators with pairwise different types
can be proved completely analogously.  Indeed, if we assume that
we can write the zero operator as a linear combination of
elementary almost natural operators of different types, then we
take the minimal
cardinality of all sets occurring in the types of the operators and
consider those which have one set of this cardinality in their type,
and so on, exactly as above. At the end we choose $r+1$ disjoint open
sets $U_i\subset M$ and sections supported in the appropriate $U_i$.
The same arguments as above then apply.
\qed\enddemo

\proclaim{\nmb.{5.5}. Corollary}
Let $D$ be an almost natural operator on a manifold $M$, $U$ an open
subset of $M$ and $I$ any subset of $\{ 1,\dots ,k\}$. If the
operator $(D|_U)^I$ associated to the restriction of $D$ to $U$ is
local then $D^I$ is local.
\endproclaim
\demo{Proof}
This is clear since by \nmb!{5.3} and \nmb!{5.4} the locality
properties of an almost natural operator are independent of the
manifold on which the operator is defined.
\qed\enddemo

\subhead \nmb.{5.6}\endsubhead
Now we can pass to the proof of our main result, the theorem
\nmb!{1.9}:
\proclaim{Theorem}
Let $E_1,\dots ,E_k,E$ be natural vector bundles defined on the
category $\Mf_m^+$, $m\geq 2$, and let $M$ be a connected oriented
smooth manifold of dimension $m$. Then any separately continuous
$k$--linear operator $D:\Ga_c(E_1M)\x\dots \x\Ga_c(E_kM)\to \Ga(EM)$
which commutes with Lie derivatives is an almost natural operator. In
particular any such operator is automatically jointly continuous and
the space of such operators is always finite dimensional and independent
of the manifold $M$.
\endproclaim
\demo{Proof}
We proceed by induction on $k$. For $k=1$ the theorem has been proved
in section \nmb!{4}. So let us assume that $k\geq 2$ and that the
theorem has been proved for all $\ell<k$. Using the induction
hypothesis, the space of absolutely invariant elements in the
natural presheaf
$U\mapsto L(\Ga_c(E_{i_1}U),\dots ,\Ga_c(E_{i_\ell}U);\Ga(EU))$
is finite dimensional and each of them is determined by any restriction
to an open subset. Thus, the presheaf satisfies condition (ii) of
\nmb!{2.1}.(3). Furthermore, exactly as in the proof of \nmb!{4.6} we
verify condition (i) of \nmb!{2.1}.(3). Hence the
natural presheaf
$U\mapsto L(\Ga_c(E_{i_1}U),\dots ,\Ga_c(E_{i_\ell}U);\Ga(EU))$  is
admissible for all $\ell<k$. So we may apply theorem \nmb!{2.9}
to the continuous linear operator
$D^1:\Ga_c(E_1M)\to L(\Ga_c(E_2M),\dots ,\Ga_c(E_kM);\Ga(EM))$. Thus
we see that $D^1=\tilde D^1+\la ^1$, where $\tilde D^1$ is a local
operator and $\la ^1$ is a continuous linear functional with values
in the space of absolutely invariant elements of
$L(\Ga_c(E_2M),\dots ,\Ga_c(E_kM);\Ga(EM))$, which are by definition
exactly the operators which commute with Lie derivatives. By
induction this is the finite dimensional space of almost natural
operators, so using proposition \nmb!{3.3} we see that $\la^1$ gives
an almost natural operator.

So let us assume that $D$ has the property that $D^1$ is local. Then
consider the associated continuous linear operator $D^2$. Again by
\nmb!{2.9} this operator splits as $D^2=\tilde D^2+\la^2$, where
$\tilde D^2$ is local and $\la^2$ again gives rise to some almost
natural operator. Next assume that $U\neq \emptyset$ is open and that
$s_2\in \Ga_c(E_2M)$ has support contained in $M\setminus \bar U$.
Then by lemma
\nmb!{2.7}, $r^M_U(\la^2 (s_2))=r^M_U(D^2(s_2))$. Then suppose that we
take a section $s_1\in \Ga_c(E_1M)$ which has support contained in
$U$ and vanishes on some open
subset $V$ of $U$ and sections $s_i\in \Ga_c(E_iM)$ with support
contained in $V$ for $i=3,\dots ,k$. Then by construction $\supp
(s_i)\subset M\setminus \supp(s_1)$ for all $i>1$, so
$D(s_1,\dots ,s_k)=D^2(s_2)(s_1,s_3,\dots ,s_k)$ restricts to zero on
$M\setminus \supp(s_1)$ and thus we see that $\la^2
(s_2)(s_1,s_3,\dots ,s_k)$ restricts to zero on $V$. Thus
$\la^2(s_2)$ restricts on $U$ to an almost natural operator which has
the property that the associated operator of type 1 is local, so by
\nmb!{5.5} $\la^2(s_2)$ has this property. Since any section $s_2$
can be written as the sum of two sections having supports in some open
subsets $U\neq M$ this holds for any $s_2$. Consequently the
$k$--linear operator associated to $\la^2$ also has the property that
its associated operator of type 1 is local and thus the same holds
for the $k$--linear operator associated to $\tilde D_2$.

Together we see that subtracting almost natural operators we come
from the original operator $D$ to an operator $\hat D$ such that
the associated operators $\hat D^1$ and $\hat D^2$ are local.

Iterating this procedure we see that subtracting suitable almost
natural operators from $D$ we arrive at an operator which has the
property that all associated linear operators are local. So let us
assume that $D$ itself already has this property.

Then we consider the first associated bilinear operator $D^{1,2}$.
By theorem \nmb!{2.9} this operator splits as
$D^{1,2}=\tilde D^{1,2}+\la^{1,2}$ where $\tilde D^{1,2}$ is local
and $\la ^{1,2}$ is weakly local. Using the induction hypothesis and
\nmb!{3.3} we see that $\la^{1,2}$ gives again rise to an almost
natural operator. Similar arguments as above show that the
$k$--linear operator corresponding to $\tilde D^{1,2}$ has the
property that all associated linear operators and the associated
bilinear operator of type $(1,2)$ are local. Let us again write $D$
for this operator.

Then we consider the associated operator $D^{1,3}$. By \nmb!{2.9}
this splits as $D^{1,3}=\tilde D^{1,3}+\la^{1,3}$, with
$\tilde D^{1,3}$ local and $\la^{1,3}$ weakly local, and again
$\la^{1,3}$ gives rise to an almost natural operator. Now
consider the $k$--linear operator $\tilde D$ associated to
$\tilde D^{1,3}$. We claim that its associated operator of type
$(1,2)$ is local. Let us first assume that we have a section $s_1$
which vanishes on some open subset $U$ of $M$ an sections $s_i$ with
support in $U$ for $i\geq 3$. Then in particular $s_1$ and $s_3$ have
disjoint supports and hence $\la^{1,3}(s_1,s_3)=0$ by weak locality.
Consequently for such sections $\tilde D$ coincides with $D$ and thus
vanishes on $U$. Moreover similarly as before one shows that
$\la^{1,3}$ has values in the space of almost natural operators which
have the property that the associated linear operator corresponding
to $s_2$ is local. Thus if we have a section $s_2$ which vanishes on
some open subset $U$ while the sections $s_i$ have support in $U$ for
all $i\geq 3$ then $\la^{1,3}(s_1,s_3)(s_2,s_4,\dots ,s_k)$ restricts
to zero on $U$. But since $D^{1,2}$ is local the same is true for
$D^{1,3}(s_1,s_3)(s_2,s_4,\dots ,s_k)$ and thus also
$\tilde D^{1,3}(s_1,s_3)(s_2,s_4,\dots ,s_k)$ restricts to zero on $U$. So
again no locality properties are lost.

Iterating this procedure we see that subtracting suitable almost
natural operators we arrive at an operator such that all associated
bilinear (and thus all associated linear operators) are local.
Next we consider the associated trilinear operators in some
succession, then the 4--linear ones and so on until we arrive at an
operator which has the property that all associated $(k-1)$--linear
operators are local. To this operator we can now directly apply
theorem \nmb!{2.9} to split it into a local and a weakly local part
which by \nmb!{3.3} is again an almost natural operator.

Finally from \nmb!{1.6} it is clear that the space of almost natural
$k$--linear operators is again finite dimensional.
\qed\enddemo

\head\totoc\nmb0{6}. Examples and Applications \endhead

In order to illustrate the strength of our results, we shall
discuss several concrete geometrical problems. In particular we
shall give a completely algebraic characterization of several
well known brackets, the Lie bracket, the Schouten bracket, the
Schouten--Nijenhuis bracket and a bracket closely related to the
Fr\"olicher--Nijenhuis bracket.

\subheading{\nmb.{6.1}. Remark}
In fact the only point where we need the continuity assumption
on the operators is the
description of weakly local multilinear functionals which commute
with Lie derivatives. Hence using assumptions which make sure that
only linear functionals can occur we can obtain automatic continuity
results. We carry this out only for the bilinear case.

\proclaim{\nmb.{6.2}. Corollary}
Let $E_1,E_2$ and $F$ be natural vector bundles defined on
$\Cal Mf_m^+$ such that $F$ has no absolutely invariant sections. Let
$M$ be a smooth oriented manifold of dimension $m\geq 2$ and let
$D:\Ga_c(E_1M)\x \Ga_c(E_2M)\to \Ga(FM)$ be a (not necessarily
continuous) bilinear operator which commutes with Lie derivatives. Then
$D$ is an almost natural operator and thus in particular jointly
continuous.
\endproclaim
\demo{Proof}
Since we proved the description of the linear weakly local
functionals commuting with Lie derivatives without continuity assumptions the
proof of \nmb!{4.6} shows also that the presheaf
$\Cal F(U):=\tilde L(\Ga_c(E_iU),\Ga (FU))$ of not necessarily
continuous linear operators is admissible. Now as in the proof of the
main theorem subtracting almost natural operators from $D$ we arrive
at an operator such that both associated linear operators are local.
Applying to this operator theorem \nmb!{2.9} we see that it must be
local since $F$ has no absolutely invariant sections.
\qed\enddemo

\subhead\nmb.{6.3}.  The operators on exterior forms\endsubhead Our first
application will give a complete description of all $k$--linear
separately continuous operators $D\:\Ga_c(\La^{p_1}T^*M)\x \dots\x
\Ga_c(\La^{p_k}T^*M)\to \Ga(\La^qT^*M)$ commuting with the Lie
derivatives, where $M$ is an arbitrary $m$--dimensional oriented
connected manifold, $m\ge 2$.

According to theorem \nmb!{1.9}, $D$ must be an almost natural
operator. Thus all operators in question are linearly generated
by the possible elementary almost natural operators. So let us
focus on the multilinear functionals $\la$ and local natural operators
$\tilde D$ which may appear in $D$. First assume $q=0$ and assume
$D=\la$ is
a multilinear functional vanishing on Lie derivatives. Then $\la$
is defined by means of a local natural operator $\hat D\:
\Ga_c(\La^{p_1}T^*M)\x \dots\x \Ga_c(\La^{p_{k -1}}T^*M)\to
\Ga_c(\La^{p_k} TM\otimes\La^mT^*M)$. However, $\La^{p_k}\Bbb
R^m\otimes\La^m\Bbb R^{m*}\simeq \La^{m-p_k}\Bbb R^{m*}$ as
$GL(m,\Bbb R)$--spaces and so the corresponding natural vector
bundles are naturally equivalent. Moreover, under this
identification, the canonical pairing with values in the
$m$--forms is just the wedge product.

What remains is to classify all local natural operators $D\:
\Ga_c(\La^{p_1}T^*M)\x \dots\x \Ga_c(\La^{p_{\ell }}T^*M)\to
\Ga_c(\La^qT^*M)$, $0\le q\le m$. This can be done very easily
by the general technique described in \nmb!{1.5}. Let us sketch
briefly how to do it. At the level of standard fibers of the jet
spaces, we can write the operator in question as a linear
combination of terms $\ph_{i_1\dots i_{p_1},\al}\dots
\ps_{j_1\dots j_{p_\ell },\be}$ denoting tensor product of
derivatives of the tensors. The indices indicated by $\al,
\dots, \be$ are symmetric, we have to alternate all indices at
the end, and there are no indices to contract over, which implies
that each individual term in the tensor product can be differentiated at
most once. Thus,
up to scalar multiples, the only way to get a local natural operator is
to choose $n=q-(p_1+\dots+ p_{\ell })\ge 0$ indices $1\le j_1<\dots
j_n\le \ell $ and to define
$$
D^{p_1,\dots,p_\ell }_{j_{1},\dots, j_n}(s_1,\dots,s_\ell ) :=
ds_{j_1}\wedge\dots\wedge ds_{j_n}\wedge
s_{i_1}\wedge\dots\wedge s_{i_{\ell -n}}
$$
where $i_1, \dots, i_{\ell -n}$ are the remaining indices among
$\{1,\dots, \ell \}$.

In words, the operators in question are built from the wedge
products, exterior differentials and integration of $m$--forms.
The classification in each concrete situation is a matter of
simple combinatorics. Let us formulate two corollaries.

\proclaim{\nmb.{6.4}. Corollary}
Let  $M$ be an arbitrary $m$--dimensional oriented
smooth manifold, $m\ge 2$. All bilinear (not necessarily
continuous) operators $D\:\Ga_c(\La^{p_1}T^*M)\x \Ga_c(\La^{p_2}T^*M)
\to \Ga(\La^qT^*M)$, $q>0$, $p_1<m$, $p_2<m$, commuting with the Lie
derivatives are local natural operators (thus obtained by means of
wedge product and exterior differential).\endproclaim

If we drop one of the assumptions that $p_1,p_2<m$ and $q>0$ then
there are operators like
$$\gather
D\:\Ga_c(\La^mT^*M)\x \Ga_c(\La^pT^*M)\to
\Ga(\La^{p+1}T^*M),\quad
(\om, \ph)\mapsto (\int_M \om)d\ph\\
D\:\Ga_c(\La^{m-p}T^*M)\x \Ga_c(\La^pT^*M)\to C^\infty(M),\quad
(\om,\ph)\mapsto \int_M(\om\wedge\ph)
.\endgather$$

\proclaim{\nmb.{6.5}. Corollary} The space of all separately
continuous $k$--linear operators
$$C^\infty_c(M)\x\dots\x C^\infty_c(M)\to C^\infty(M)$$
is linearly generated by the elementary almost natural operators
obtained by choosing pairwise disjoint $m$-tuples of indices
$I_j=\{i^j_1<\dots<i^j_m\}$ among $\{1,\dots,k\}$, considering the
$m$--forms $df_{i^j_1}\wedge\dots\wedge df_{i^j_m}$, distributing
some of the other functions to the $m$--forms in order to
integrate them, and to multiply the result of the integrations
with the remaining functions.

In particular, all symmetric separately continuous $k$--linear
operators are the scalar
multiples of the pointwise multiplication.
\endproclaim

This corollary is closely related to one of the results in \cite{de
Wilde--Lecomte}, where the last assertion of the
corollary is proved without any continuity assumptions.

\subheading{\nmb.{6.6} }
A more difficult task is to get similar complete descriptions for
operations with arguments involving some contravariant indices.
We shall restrict ourselves to some bilinear operations.
First we reprove the completely algebraic characterization of the
Lie bracket, cf\. \cite{de~Wilde--Lecomte}.
\proclaim{Theorem}
Let $M$ be an oriented smooth manifold of dimension $\geq 2$,
$D:\Cal X_c(M)\x\Cal X_c(M)\to \Cal X(M)$ a (not necessarily
continuous) bilinear operator which maps compactly supported vector
fields to vector fields and commutes with Lie derivatives. Then $D$
is a scalar multiple of the Lie bracket.
\endproclaim
\demo{Proof}
Since $TM$ has no absolutely invariant sections,
$D$ must be an almost natural operator
according to \nmb!{6.2}. Moreover
$T^*M\otimes \La^mT^*M$ has no absolutely invariant sections since
the elements $\la \cdot \text{Id}$ from the center of $GL(m,\Bbb
R)$ act by multiplication with a nonzero power of $\la$, and
so $D$ must in fact be
a local natural operator. Now assuming locality it is easy to prove
uniqueness of the Lie bracket by means of the technique indicated
in \nmb!{1.5} (c\.f\. \cite{\v Cap},
\cite{Krupka--Mikol\'a\v sov\'a} and \cite{van~Strien}).
\qed\enddemo

\subhead \nmb.{6.7}. The Schouten bracket\endsubhead
The Lie bracket is a special case of the so called Schouten
bracket $D\:\Ga_c(S^kTM)\x
\Ga_c(S^\ell TM)\to \Ga(S^{k+\ell-1}TM)$, $k>0$, $l>0$,
(or symmetric Schouten concomitant) which can be defined as
the restriction of the canonical Poisson bracket on the
symplectic manifold $T^*M$ to the fiberwise polynomial functions
on $T^*M$ identified with sections of $S^kTM$ and $S^\ell TM$.

Since the natural bundle $S^kTM$, $k>0$, corresponds to a
non-trivial irreducible representation of $GL(m,\Bbb R)$, there
are no absolutely invariant sections in the target of the
operations in question. Thus \nmb!{6.2} implies that the operator
must be a linear combination of elementary almost natural
operators. But since there are no absolutely invariant sections
in $S^kT^*M\otimes \La^mT^*M$, the operator must be local. Then
the general procedure for the classification shows that there is
a unique local operator, up to scalar multiples. Thus we have
proved

\proclaim{Theorem} Let $M$ be an $m$--dimensional oriented
smooth manifold, $m\ge 2$. Each bilinear operator
$D\:\Ga_c(S^kTM)\x \Ga_c(S^\ell TM)\to \Ga(S^{k+\ell-1}TM)$, $k>0$, $l>0$,
which commutes with Lie derivatives is a scalar multiple of the
Schouten bracket.
\endproclaim

\subhead \nmb.{6.8}. The Schouten--Nijenhuis bracket\endsubhead
Another generalization of the Lie bracket is the bracket
$D\:\Ga_c(\La^pTM)\x \Ga_c(\La^qTM)\to \Ga_c(\La^{p+q-1}TM)$
defined by the formula
$$\multline
D(X_1\wedge\dots\wedge X_p, Y_1\wedge\dots\wedge Y_q)=
\\
= \sum_{i,j} (-1)^{i+j}[X_i,Y_j]\wedge X_1\wedge\dots
{}^{\wedge_i}\ldots\wedge
X_p\wedge Y_1\wedge\dots{}^{\wedge_j}\ldots\wedge Y_q.
\endmultline$$
This operation is called Schouten--Nijenhuis bracket and defines
the structure of a graded Lie algebra on the space $\La TM$.

Exactly as in 5.10, both  the target bundle and $\La^qT^*M\otimes
\La^mT^*M$ admit no absolute invariant sections and so all
operators of the given type which commute with Lie derivatives
(without any continuity assumption) must be local. But this
implies, cf\. \cite{Michor 1},

\proclaim{Theorem}
Let $M$ be an $m$--dimensional oriented
smooth manifold, $m\ge 2$. Each bilinear operator
$D\:\Ga_c(\La^pTM)\x \Ga_c(\La^q TM)\to \Ga(^{p+q-1}TM)$, $p>0$,
$q>0$,
which commutes with Lie derivatives is a scalar multiple of the
Schouten--Nijenhuis bracket.
\endproclaim

\subheading{\nmb.{6.9}. The Fr\"olicher--Nijenhuis bracket}
Consider the first order natural
vector bundles $\La^kT^*\otimes T$. Sections of values of such a
natural bundle are
usually called vector valued $k$--forms, and the space of these
sections over $M$ is denoted by $\Om^k(M,TM)$. We shall give a
completely algebraic characterization of a bracket
which is closely related to the Fr\"olicher--Nijenhuis bracket.

For any $M$ the bundle
$\La T^*M\otimes TM=\oplus_k\La^kT^*M\otimes TM$ is a bundle of graded
modules over the bundle of graded commutative algebras $\La T^*M$,
where the action is induced by the wedge product. Moreover there is
an absolutely invariant section $\Bbb I$ of $T^*M\otimes TM$, which
corresponds to the identity map on $TM$. By acting on this element we
can view $\La^{k-1}T^*M$ as a subbundle of $\La^kT^*M\otimes TM$.
This subbundle is actually a direct summand, and a projection
$\La^kT^*M\otimes TM\to \La^{k-1}T^*M$ is induced by mapping
$\om\otimes X$ to $\tfrac{(-1)^{k-1}}{m-k+1}i_X\om$, where $i_X$
denotes the usual insertion operator. The kernel of this projection
is again a natural vector bundle $C^k$, whose sections are called trace
free vector valued $k$--forms. We shall write briefly
$\Cal C^k(M)$ for the space
of these sections. (The decomposition of $\La^kT^*M\otimes TM$
constructed above corresponds exactly to the decomposition of the
representation $\La ^k\Bbb R^{m*}\otimes \Bbb R^m$ of $GL(m,\Bbb R)$
into irreducible representations.)

There is a bilinear natural operator
$[\quad,\quad ]:\La^kT^*\otimes T\x \La^\ell T^*\otimes T\to
\La^{k+\ell} T^*\otimes T$ called the Fr\"olicher--Nijenhuis bracket,
which defines a graded Lie algebra structure on the space of vector
valued differential forms and plays an important role in the theory
of generalized connections (c\.f\. \cite{Fr\"olicher--Nijenhuis} and
\cite{Michor 2}). Using the projections onto trace free
vector valued forms constructed above we can compress this bracket to
a bilinear natural operator
$[\quad,\quad]^c:\Cal C^k(M)\x \Cal C^\ell(M)\to C^{k+\ell}(M)$ for
any smooth manifold $M$. It turns out that this is again a graded Lie
bracket (c\.f\. \cite{Michor--Schicketanz}).

\proclaim{\nmb.{6.10}. Theorem}
Let $M$ be an oriented smooth manifold of dimension $m\geq 2$,
$D:\Cal C_c^k(M)\x\Cal C_c^\ell (M)\to \Cal C^{k+\ell}(M)$ a not
necessarily continuous bilinear operator which commutes with Lie
derivatives. Then $D$ is a scalar multiple of the compression of the
Fr\"olicher--Nijenhuis bracket.
\endproclaim
\demo{Proof}
First note that $\Cal C^{k+\ell}$ has no nonzero absolutely invariant
sections since it corresponds to an irreducible representation of
$GL(m,\Bbb R)$. Thus by \nmb!{6.2} $D$ must be an almost natural
operator. Moreover for $k<m$ there are no absolutely invariant
sections of $C^{k*}M\otimes \La^mT^*M$, since scalar multiples of the
identity act nontrivially. But $\Cal C^m(M)=0$ since $\La^m\Bbb
R^{m*}\otimes\Bbb R^m\simeq \La^{m-1}\Bbb R^{m*}$ (c\.f\.
\cite{Michor--Schicketanz}). Thus $D$ must be a local operator.
Now using the projections on the left hand side and the inclusion
into all vector valued differential forms on the right hand side, $D$
induces a local natural operator
$\tilde D:\Om^k(M,TM)\x \Om^\ell(M,TM)\to \Om^{k+\ell}(M,TM)$.
But these operators have been completely classified in
\cite{Kol\'a\v r--Michor} and \cite{\v Cap}. Looking at the
list of possible operators and taking into account that the values
must be trace free and the operator depends only on the trace free
parts of the arguments one sees that in fact it must be a scalar
multiple of the Fr\"olicher--Nijenhuis bracket, so the result
follows.
\qed\enddemo

%

\enddocument